\magnification=\magstep1
\overfullrule=0pt
\baselineskip=20truept
\def\a{\alpha}
\def\b{\beta}

\def\eip{e^{i\phi}}
\def\meip{e^{-i\phi}}
\def\eit{e^{i\theta}}
\def\meit{e^{-i\theta}}
\def\G{\Gamma}
\def\K{{\cal K}}
\def\bg#1{\Big({#1}\Big)}
\def\d{\delta}
\def\g{\gamma}
\def\o{\over}

\def\itt{\int\limits^\infty_{-\infty}}

\def\ty{\infty}
\def\ta{\theta}
\def\sq#1{\sqrt{#1}}

\def\lm{\lambda}
\def\sn{\sum\limits^\infty_{n=0}}
\def\nk{\sum\limits^n_{k=0}}
\def\tk{\sum\limits^\infty_{k=0}}
\def\ep{\epsilon}

\def\om{\omega}

%\def\s{\sigma}
%\bigskip
\vskip 2in
\centerline{\bf On a general $q$-Fourier transformation}
\centerline{\bf with nonsymmetric kernels\footnote {\dag}{\rm This
work was supported in part by NSF grant \#DMS 93005224 and
the Natural Sciences and Engineering
Research Council of Canada grant \#A6197.}}
\bigskip
\centerline{Richard A. Askey\footnote*{Dept. of Mathematics,
University of
Wisconsin, Madison, Wisconsin 53706}, Mizan
Rahman\footnote{**}{Dept. of
Mathematics and
Statistics, Carleton University, 1125 Colonel By Drive, Ottawa,
Ontario, Canada K1S 5B6}
and Serge\u \i\ K. Suslov\footnote{***}{Russian Research Centre
``Kurchatov Institute'', Moscow 123182, Russia.}}
\bigskip
\centerline{\bf Abstract}
\bigskip
Wiener used the Poisson kernel for the Hermite polynomials to deal
with the classical Fourier transform. Askey, Atakishiyev and Suslov
used this approach to obtain a $q$-Fourier transform by using the
continuous $q$-Hermite polynomials. Rahman and Suslov extended this
result by taking the Askey--Wilson polynomials, considered to be
the
most general continuous classical orthogonal polynomials. The
theory of $q$-Fourier transformation is further extended here by
considering a nonsymmetric version of the Poisson kernel with
Askey--Wilson polynomials. This approach enables us to obtain some
new results, for example, the complex and real orthogonalities of
these kernels.
\vfil\eject
\noindent
{\bf 1. Introduction}.\quad Two of the most frequently used
formulas in the area of integral transforms are the classical
Fourier
transform
$$f(x) = {1\o \sq{2\pi}} \itt e^{ixy} g(y) \, dy\ \colon =
F[g](x),\leqno (1.1)
$$
and the corresponding formal inversion formula
$$g(y) = {1\o \sq{2\pi}} \itt e^{-ixy} f(x)\, dx\ \colon = F^{-1}
[f]
(y).\leqno (1.2)
$$
There are many excellent references on Fourier transforms (see, for
example, Titchmarsh's classic book [20]), but, for our purposes the
most important reference is Wiener [21], who looked at the Fourier
transform from the point of view of classical orthogonal
polynomials on $R$, in particular, the Hermite polynomials
$H_n(x)$.
\medskip
An elementary result in non-relativistic Quantum Mechanics is that
the normalized wavefunction of the Hamiltonian operator for the
harmonic oscillator is the Hermite function
$$\Psi_n(x) = \bg{2^nn!\sq{\pi}}^{-1/2} H_n (x)\, e^{-x^2/2},\leqno
(1.3)$$
see, for example, Landau and Liftschitz [13, p.70]. The bilinear
generating function (the Poisson kernel) for these functions is
$$\leqalignno{
K_t(x,y) &= \sn t^n \, \Psi_n (x) \, \Psi_n(y)&(1.4)\cr
&= [\pi (1-t^2)]^{-1/2} \exp \Big[{4xyt - (x^2 + y^2) (1+t^2) \o
2(1-t^2)}\Big],\cr
}
$$
which is better known as Mehler's formula, see [19]--[21].
\medskip
We make two important observations. First, the kernel $K_t(x,y)$
as given by the function on the right of the second equality in
(1.4) has the
property
$$\lim_{t\to i} K_t (x, y) = K_i (x,y) = {e^{ixy}\o
\sq{2\pi}},\leqno (1.5)$$
which, of course, is the kernel of the integral transform (1.1). In
fact, $K_t (x,y)$ is analytic everywhere on the unit circle in the
complex $t$-plane except at $t = \pm 1$. We shall set
$$\leqalignno{
\K_\tau (x,y) &= K_t (x, y)\Big|_{t = e^{i\tau}}&(1.6)\cr
&= {e^{i(\pi/2 - \tau)/2}\o \sq{2\pi\sin\tau}} \exp \Big[i {2xy -
(x^2+ y^2) \cos\tau\o 2\sin\tau}\Big],\cr
}
$$
$0 < \tau < \pi$, and define a generalized Fourier transform by the
formula
$$f(x) \colon = \itt \K_\tau (x,y) \, g(y)\, dy.\leqno (1.7)$$
The corresponding formal inversion formula is
$$g(y) = \itt \K^*_\tau (x,y) \, f(x)\, dx,\leqno (1.8)$$
where $*$ indicates the complex conjugate. For a sketch of the
proof of (1.8) see [17].
\medskip
The second observation is that the wavefunctions in (1.3) are the
eigenfunctions of the Fourier operators in (1.1) and (1.7):
$$e^{in\tau} \, \Psi_n(x) = \itt \K_\tau (x,y) \, \Psi_n (y)\,
dy,\leqno
(1.9)$$
where $0 < \tau < \pi$ for (1.7) and, $\tau = \pi/2$ for (1.1).
\medskip
Our primary interest in this paper is a general $q$-analogue of
(1.9) as well as of (1.7). As was explained in [17], a
$q$-analogue of a formula contains a parameter $q$, usually complex
and $|q| < 1$, such that the limit of the new formula as $q\to 1$
is the given formula. For example, a $q$-analogue of a complex
number $a$ is $(1-q^a)/(1-q)$ since
$$a = \lim_{q\to 1} {1-q^a\o 1-q},$$
the branch of $q^a$ chosen here is the one that gives
$\lim\limits_{q\to 1} q^a = 1$.
\medskip
There is a great deal of interest these days in $q$-analogues of
important classical formulas. Sometimes the interest is nothing
more than a curiosity, but in view of some recent developments in
the generalizations of the classical harmonic oscillator problem
(see, for example, [4,5], [8,9] and [14]) the interest in a
$q$-version of the Fourier transform is much more than an academic
curiosity. The question appears in a very natural way in the
$q$-oscillator problem.
\medskip
As was shown in [2] it is Wiener's treatment of the Fourier
integrals that contains the key to a meaningful $q$-extension. Just
as the Hermite polynomials are associated with the wavefunctions
for the harmonic oscillator, the continuous $q$-Hermite
polynomials,
$$H_n(x|q) =\nk {(q;q)_n\o (q;q)_k (q;q)_{n-k}}\,
e^{i(n-2k)\ta},\quad x = \cos\ta,\leqno (1.10)$$
are associated with the $q$-wave functions for the $q$-harmonic
oscillator, namely,
$$\Psi_n (x|q) = \Big[ (q^{n+1};q)_\ty\Big/2\pi\Big]^{1/2}
\sq{\rho_0(x)}\, H_n(x|q),\leqno (1.11)
$$
where
$$\rho_0 (x) = 4\sqrt{1-x^2}\, \prod^\ty_{k=1} \bg{1-2 (2x^2-1) q^k
+ q^{2k}},\ 0 < q < 1. \leqno (1.12)$$
The $q$-shifted factorials in (1.10) and (1.11) are defined by
$$\eqalign{
&(a;q)_0 = 1, \quad (a;q)_n = (1-a) (1-aq)\ldots (1-aq^{n-1}),\quad
n=1,2, \ldots, \cr
&(a;q)_\ty = \lim_{n\to\ty} (a;q)_n, \quad |q|< 1.\cr} \leqno
(1.13)$$
The $q$-annihilation operator $b$ and the $q$-creation operator
$b^+$ that satisfy the commutation rule
$$bb^+ - q^{-1} b^+b = 1\leqno (1.14)$$
were introduced explicitly in [8]. They act on the $q$-wavefunction
(1.11) in the following manner:
$$\eqalign{
&b \Psi_n (x|q) = \bg{1-q^{-n}\o 1-q^{-1}}^{1/2} \Psi_{n-1}
(x|q),\cr
&b^+ \Psi_n(x|q) = \bg{1-q^{-n-1}\o 1-q^{-1}}^{1/2} \Psi_{n+1}
(x|q),\cr}\leqno (1.15)
$$
see [2] and [8]. One also has the orthogonality property
$$\int^1_{-1} \Psi_m (x|q) \, \Psi_n(x|q)\, dx =
\delta_{m,n},\leqno
(1.16)$$
see, for example, [3] and [10]. The important difference between
this formula and the corresponding formula for the Hermite
functions (1.3), namely,
$$\itt \Psi_m (x) \, \Psi_n (x)\, dx = \delta_{m,n},\leqno (1.17)$$
is that one is over the finite interval $(-1, 1)$ while the other
is
over the whole real line. So the use of (1.11) as an analogue of
(1.3) has the advantage of orthogonality over a finite interval.
This automatically ensures completeness of the system $\{\Psi_n
(x|q)\}^\ty_{n=0}$, see [19, Thm, 3.1.5]. Furthermore, we need not
be tied to the $q$-analogues of the Hermite polynomials only, but
can consider other generalizations of these polynomials, e.g., the
$q$-ultraspherical [3], the continuous $q$-Jacobi [6], [15] or even
the most general, the Askey--Wilson, polynomials [6], all of which
are orthogonal on $(-1, 1)$. Following this direction, the ideas
presented in  [2] were extended in [17] to the case of
the Askey--Wilson polynomials using the explicit formula for the
Poisson kernel derived in [18]. As in [2] this led to a singular
integral transformation, thus extending the results corresponding
to the continuous $q$-Hermite polynomials.
\medskip
The purpose of this paper is to do even more: to apply the ideas of
Wiener's treatment of the Fourier integrals to the case of
a nonsymmetric extension of the Poisson kernel for the
Askey--Wilson
polynomials derived here. This paper is organized as follows. In
section 2 we introduce the Askey--Wilson polynomials and then state
our result for the corresponding Poisson kernel in section 3. In
sections 4--6 we derive an explicit formula for this kernel which
reveals the structure of the poles. The orthogonality property of
these kernels is then proved in sections 7--12. Finally, we
introduce a new integral transformation and prove its inversion
formula in section 13. We close  this paper by displaying some
special cases of the kernel and a continuous orthogonality relation
for the simplest of them, in section 14.
\bigskip
\noindent
{\bf 2. The Askey--Wilson polynomials}.\quad We shall assume that
$a, b, c, d$ and $q$ are real parameters such that $0 < q < 1$ and
$\max (|a|, |b|, |c|, |d|)<1$. The Askey--Wilson polynomials [6]
are
defined by
$$\leqalignno{
p_n(x) &= p_n (x; a, b, c, d)&(2.1)\cr
&= \, _4\phi_3\left[\matrix{
q^{-n}, \, abcdq^{n-1}, \, ae^{i\ta}, \, ae^{-i\ta}\cr
ab, \, ac, \, ad\cr} ;q,q\right],\cr
}
$$
$n = 0, 1, \ldots, x = \cos\ta$, $0 \le \ta \le \pi$. This
definition differs slightly from the one given in [6], but, for our
purposes, this is the more convenient one. The symbol on the right
side is the $r = 3$ case of the basic hypergeometric series
$_{r+1}\phi_r$ defined by
$$\leqalignno{
&_{r+1}\phi_r\left[\matrix{
a_1, a_2, \ldots, a_{r+1}\cr
b_1, b_2, \ldots, b_r\cr} ; \,q,\, z\right]&(2.2)\cr
&\quad= \sum^\ty_{k=0} {(a_1, a_2, \ldots, a_{r+1};q)_k\o
(q, b_1, \ldots, b_r;q)_k}\, z^k,\cr
}
$$
with
$$(a_1, a_2, \ldots, a_m;q)_k = \prod^m_{j=1} (a_j;q)_k, \quad
\hbox{see [10]}.\leqno(2.3)$$
The orthogonality property satisfied by $p_n (x; a, b, c, d)$ is
$$\int^1_{-1} p_n (x) \, p_m (x) \, \rho (x)\, dx =
\delta_{m,n}/h_n,\leqno
(2.4)
$$
$$
\leqalignno{
&h_0 = {(q, ab, ac, ad, bc, bd, cd;q)_\ty\o
2\pi (abcd;q)_\ty}, &(2.5)\cr
&h_n = h_0\, {(1-abcdq^{2n-1}) (abcdq^{-1}, ab, ac, ad;q)_n\o
(1-abcdq^{-1}) (q, cd, bd, bc;q)_n}\, a^{-2n},&(2.6)\cr
&\rho(x) = \rho (x; a,b, c, d) = {h(x; 1, -1, q^{1/2}, -q^{1/2})\o
h(x; a, b, c, d)}\, (1-x^2)^{-1/2},&(2.7)\cr
\noalign{\hbox{and}}
& h(x; a_1, a_2, \ldots, a_r) = \prod^r_{j=1} h(x; a_j),&(2.8)\cr
&h(x;a) = \prod^\ty_{n=0} (1-2axq^n+ a^2 q^{2n}).&(2.9)\cr
}
$$
{\bf 3. The Poisson kernel and a nonsymmetric extension}.\quad The
Poisson kernel for the Askey--Wilson polynomials is defined by
$$\leqalignno{
P_t (x,y) &= h_0\, K_t (x,y),&(3.1)\cr
\noalign{\hbox{with}}
K_t(x,y) =  \sn &{(1-abcdq^{2n-1})(abcdq^{-1}, ab, ac, ad;q)_n\o
(1-abcdq^{-1}) (q, cd, bd, bc;q)_n}\, a^{-2n} t^n&(3.2)\cr
&\cdot p_n (x; a, b, c, d)\, p_n (y; a, b, c, d).\cr
}
$$
The series on the right side of (3.2) converges for $|t|<1$ and
$|x| \le 1$, $|y|\le 1$. In [11], Gasper and Rahman found an
explicit
representation of $K_t(x,y)$ in the special case $ad = bc$, and
proved its positivity when $a = q^{\a/2 + 1/4}$,  $b = aq^{1/2}$,
$c = -q^{\b/2 + 1/4}$, $d = cq^{1/2}$, $\a, \b > -1$ (this
corresponds to Askey and Wilson's continuous $q$-Jacobi polynomials
[6]). Rahman and Verma [18] found an expression for the most
general kernel (without the condition $ad = bc$) which was
reproduced in [17] since the original expression in [18] has many
misprints.
\medskip
In this paper we shall consider a nonsymmetric extension of the
bilinear generating function (3.2) by replacing $p_n (y; a, b,
c,d)$ by $p_n (y; \a, \b, \gamma, \delta)$ such that
$$\a\g = ac,\quad \b\delta = bd,\quad  \max(|\a|, |\b|, |\g|,
|\delta|)
< 1.\leqno (3.3)$$
For ease of notation we shall use the labels $\lm$ and $\mu$ to
refer to the parameter--quartets $(a, b, c, d)$ and $(\a, \b,
\gamma, \delta)$, respectively. Thus $h^\lm_0$ means the same as
$h_0$ in (2.5), but
$$h_0^\mu\colon = {(q, \a\b, \a\g, \a\delta, \b\g, \b\delta,
\g\delta;q)_\ty\o
2\pi (\a\b\g\delta;q)_\ty}, \leqno (3.4)$$
with similar meanings for $h^\lm_n$ and $h^\mu_n$. Also, we shall
use the notations
$$\eqalign{
&p^\lm_n (x) = p_n (x; a, b, c, d),\cr
&p^\mu_n (y) = p_n (y; \a, \b, \g, \delta),\cr}\leqno (3.5)$$
and
$$
\ep = (abcd)^{1/2}.\leqno (3.6)$$
The extension of (3.1) and (3.2) that we have in mind is
$$\leqalignno{
&P^{\lm, \mu}_t (x,y) = h_0^\lm \, K^{\lm,\mu}_t (x,y),&(3.7)\cr
\noalign{\hbox{where}}
&K_t^{\lm, \mu} (x,y) = \left(h_0^\lm\right)^{-1}
\sn h^\lm_n\, t^n\, p^\lm_n (x)\, p^\mu_n (y).&(3.8)\cr
}
$$
We shall prove in the following sections that
$$K_t^{\lm,\mu} (x, y) = K^{(1)}_t (x,y) + K_t^{(2)} (x,y) +
K_t^{(3)} (x,y),\leqno (3.9)
$$
where
$$\leqalignno{
&K_t^{(1)} (x,y)&(3.10)\cr
&= (1-t^2)\, {(-qt\ep, \a\b c/b, \a\b d/b, ae^{i\ta}, \ep^2
e^{-i\ta}/b;q)_\ty\o
(-t/\ep, ac, ad, \a\b e^{i\ta}/b, \a\b cd e^{-i\ta}/b;q)_\ty}\cr
&\cdot \tk {(\ep,\, \ep q^{1/2},\, -\ep q^{1/2},\, -\ep/ad;q)_k \o
(q,\, bc,\, -qt\ep,\, -q\ep/t;q)_k}\,q^k\cr
&\cdot \sum^k_{\ell = 0} {(q^{-k}, -\ep, ad, \a e^{i\phi}, \a
e^{-i\phi}, \a\b e^{i\ta}/b, \a \b e^{-i\ta}/b, \a\b cd
e^{-i\ta}/b;q)_\ell\o
(q, -adq^{1-k}/\ep, \a\b, \a\delta, \a/\delta, \a\b d/b, \a\b c/b,
\ep^2 e^{-i\ta}/b;q)_\ell}\, q^\ell\cr
&\cdot \sum^\ell_{m=0} {(1-\delta q^{2m-\ell}/\a) (\d q^{-\ell}/\a,
q^{1-\ell}/\a\b, bq^{1-\ell}/\a\b c, q^{-\ell};q)_m\o
(1-\d q^{-\ell}/\a) (q, \b\d, cd, q\d/\a;q)_m}\cr
&\cdot{(\d e^{i\phi}, \d e^{-i\phi}, de^{i\ta}, de^{-i\ta};q)_m\o
(q^{1-\ell} e^{-i\phi}/\a, q^{1-\ell} e^{i\phi}/\a, bq^{1-\ell}e^{-
i\ta}/\a\b, bq^{1-\ell} e^{i\ta}/\a\b;q)_m} \, \bg{bcq\o\a\d}^m\cr
&\cdot\, _8W_7 \bg{ {\a\b cd\o b} q^{\ell-1} e^{-i\ta}; ce^{-i\ta},
dq^m e^{-i\ta}, {\a\b \o b} q^{\ell-m}e^{-i\ta}, cd q^\ell, {\a\b\o
ab}; q, ae^{i\ta}},\cr
}
$$
$$\leqalignno{
&K^{(2)}_t (x,y)&(3.11)\cr
&= {(\ep^2, -\ep/ad, t, -t\ep/ad, \a\b c/b, \a\b d/b, ae^{i\ta},
\ep^2e^{-i\ta}/b;q)_\ty\o
(-\ep, bc, t/ad, -\ep/t, ac, ad, \a\b e^{i\ta}/b, \a\b cde^{-
i\ta}/b;q)_\ty}\cr
&\cdot \tk {(-t,\, tq^{1/2},\, -tq^{1/2},\, t/ad;q)_k \o
(q,\, qt^2,\, -t\ep/ad,\, -qt/\ep;q)_k}\, q^k\cr
&\cdot \sum^\ty_{\ell = 0} {(-\ep q^{-k}/t, -\ep, ad, \a e^{i\phi},
\a e^{-i\phi}, \a\b e^{i\ta}/b, \a\b e^{-i\ta}/b, \a\b cde^{-
i\ta}/b;q)_\ell\o
(q, adq^{1-k}/t, \a\b, \a\d, \a/\d, \a\b c/b, \a\b d/b, \ep^2 e^{-
i\ta}/b;q)_\ell}\, q^\ell\cr
&\cdot\sum^\ell_{m=0} {(1-\d q^{2m-\ell}/\a) (\d q^{-\ell}/\a,
q^{1-\ell}/\a\b, bq^{1-\ell}/\a\b c, q^{-\ell};q)_m\o
(1-\d q^{-\ell}/\a) (q, \b\d, cd, q\d/\a;q)_m}\cr
&\cdot{(\d e^{i\phi}, \d e^{-i\phi}, de^{i\ta}, de^{-i\ta};q)_m\o
(q^{1-\ell} e^{-i\phi}/\a, q^{1-\ell} e^{i\phi}/\a, bq^{1-\ell}
e^{-i\ta}/\a\b, bq^{1-\ell} e^{i\ta}/\a\b;q)_m} \, \bg{bcq\o
\a\d}^m\cr
&\cdot\ _8W_7 \bg{ {\a\b cd\o b} q^{\ell-1} e^{-i\ta}; ce^{-i\ta},
dq^m e^{-i\ta}, {\a\b \o b} q^{\ell-m}e^{-i\ta}, cdq^\ell, {\a\b\o
ab}; q, ae^{i\ta}},\cr
}
$$
and
$$\leqalignno{
&K^{(3)}_t (x,y)&(3.12)\cr
&= {(\ep^2, a\meit, c\meit, d\meit, \a\b\eit/b;q)_\ty \o
(ac, bc, cd, \a\b, ad/t, bc\d t/\g;q)_\ty} \cr
&\cdot {(ct\eit, c\d t \eit/\g, bct\eip/\g, bct\meip/\g;q)_\ty \o
(e^{-2i\ta}, cte^{i(\ta + \phi)}/\g, cte^{i(\ta -
\phi)}/\g;q)_\ty}\cr
&\cdot \tk {(t, -\ep/ad, -t\ep/ad, bc\d t/\g,
cte^{i(\ta + \phi)}/\g, cte^{i(\ta - \phi)}/\g;q)_k \o
(q, qt/ad, ct\eit, c\d t\eit/\g,
bct\eip/\g, bct\meip/\g;q)_k}\, q^k \cr
&\cdot \, _4\phi_3\left[\matrix{
q^{-k},\, -t,\, tq^{1/2},\, -tq^{1/2}\cr
qt^2,\, -t\ep/ad,\, -adq^{1-k}/\ep\cr};\,q,\,q\right]\cr
&\cdot \sum^\ty_{\ell = 0}
{(ctq^k e^{i(\ta + \phi)}/\g, ctq^k e^{i(\ta - \phi)}/\g,
a\eit, c\eit, d\eit;q)_\ell \o
(ctq^k\eit, c\d tq^k \eit/\g, \a\b\eit/b, qe^{2i\ta},
q;q)_\ell}\,q^\ell\cr
&\cdot\,
_8W_7\left({bc\d \o \g}t q^{k-1}; {bc \o \b\g}tq^k, {bc \o
\a\g}tq^k,
bq^{-\ell}\meit, \d\eip, \d\meip;q,\, {\a\b \o b}q^\ell \eit
\right)\cr
&+ \ {\rm idem}\ (\ta;-\ta),\cr
}
$$
where idem $(\a;\d)$ means the same expression as the preceding
one with $\a$ and $\d$ interchanged (see [10]), with a similar
meaning for
idem $(\ta; -\ta)$. In (3.10)--(3.12), $x = \cos\ta$, $0 < \ta <
\pi$, $y = \cos\phi$, $0 < \phi < \pi$, and the $W$-series are
special cases of the very-well-poised
$_{r+1}W_r$ series defined by
$$\leqalignno{
&_{r+1} W_r \bg{a; b_1, b_2, \ldots, b_{r-2};q, z}&(3.13)\cr
&\quad=\ _{r+1}\phi_r\left[\matrix{
a, qa^{1/2}, -qa^{1/2}, b_1, b_2, \ldots, b_{r-2}\cr
a^{1/2}, -a^{1/2}, aq/b_1, aq/b_2, \ldots, aq/b_{r-2}\cr}
;q,z\right],\cr
}
$$
see [10] for further details.
\medskip
Looking at the horrendous expressions in (3.10)--(3.12) one might
feel somewhat skeptical about their usefulness but, as we shall see
later, the main purpose of these formulas is to isolate the poles
of $K_t^{\lm, \mu} (x,y)$, as a function of $t$, which are
practically all we need for our subsequent analysis.
\bigskip
\noindent
{\bf 4. The $q$-integral representations}.\quad To prove
(3.9)--(3.12)
we will use, following the method of [18], the $q$-integral
representation of the Askey--Wilson polynomials [10, Ex.~7.34]:
$$\leqalignno{
p_n(x; a,b, c, d) &= \left[A(\ta)\right]^{-1} {(bc;q)_n\o
(ad;q)_n}&(4.1)\cr
&\cdot \int^{q\meit/d}_{q\eit/d} {(du\eit, du\meit, \ep^2
u/q;q)_\ty\o
(dau/q, dbu/q, dcu/q;q)_\ty}\cr
&\cdot {(q/u;q)_n\o (\ep^2 u/q;q)_n} (adu/q)^n\, d_qu,\cr
}
$$
where
$$\leqalignno{
A(\ta) &= A(\ta; a, b, c, d)&(4.2)\cr
&= -i\, {q(1-q)\o 2d} \, (q, ab, ac, bc;q)_\ty \, h(x;d) \,
\rho (x; a, b, c, d),\cr
}
$$
and the $q$-integral is defined by
$$\eqalign{
\int^b_a f(u) \,d_q u= \int^b_0 f(u)\, d_q u - \int^a_0 f(u)\, d_q
u, &\cr
\int^a_0 f(u)\, d_q u = a(1-q) \sum^\ty_{m=0} f(aq^m)\, q^m, &\cr
}\leqno (4.3)$$
for any function $f$ such that the series on the right side of
(4.3) converge.
\medskip
By using the symmetry property of the Askey--Wilson polynomials
(2.1), see [6] and [10], and the $q$-integral representation (4.1),
we obtain
$$\leqalignno{
K_t^{\lm, \mu} (x,y) = B(\ta, \phi)
&\int^{q\meit/b}_{q\eit/b}d_qu\,
{(bu\eit, bu\meit, \ep^2 u/q;q)_\ty\o
(bau/q, bcu/q, bdu/q;q)_\ty} &(4.4)\cr
&\cdot \int^{q\meip/\g}_{q\eip/\g} d_qv \, {(\g v \eip, \g v\meip,
\ep^2
v/q;q)_\ty\o
(\g\a v/q, \g\b v/q, \g\d v/q;q)} \cr
&\cdot\, _6W_5 \left( \ep^2/q; ad, q/u, q/v;q,
bcuvt/q^2\right),\cr
}
$$
where the $\lm$ and $\mu$ parameters are related by (3.3) and, of
course, $\ep^2 = abcd = \a\b\g\d$. Also,
$$\leqalignno{
B^{-1} (\ta,\phi) = - \  & {q^2 (1-q)^2\o 4b\g} \,(q, q, ac, ad,
cd,
\a\b, \a\d, \b\d;q)_\ty&(4.5)\cr
&\cdot h(x;b)\, h(y;\g)\, \rho^\lm (x)\, \rho^\mu (y),\cr
}
$$
where $\rho^\lm (x) = \rho (x; a, b, c, d)$ and $\rho^\mu(y) = \rho
(y; \a, \b, \g, \d)$.
\medskip
Crucial to our calculations is the following formula for the
$_6W_5$ series in (4.4) which was obtained in [18]:
$$\leqalignno{
&_6W_5 \left( \ep^2/q; ad, q/u, q/v;q, bcuvt/q^2 \right)&(4.6)\cr
&= (1-t^2) \,{(-qt\ep;q)_\ty\o (-t/\ep;q)_\ty}\, \tk {(\ep, \ep
q^{1/2}, -\ep q^{1/2}, -\ep/ad;q)_k \o
(q, bc, -qt\ep, -q\ep/t;q)_k}\, q^k\cr
&\qquad\qquad \cdot \, _4\phi_3\left[\matrix{
q^{-k},\, -\ep,\, ad,\, \ep^2 uv/q^2\cr
-adq^{1-k}/\ep,\, u\ep^2/q,\, v\ep^2/q\cr};q,q\right]\cr
&+ {(\ep^2, - \ep/ad, t, -t\ep/ad;q)_\ty \o
(-\ep, bc, t/ad, -\ep/t;q)_\ty}\, \tk {(-t, tq^{1/2}, -tq^{1/2},
t/ad;q)_k\o
(q, qt^2, -t\ep/ad, -qt/\ep;q)_k}\, q^k\cr
&\qquad\qquad \cdot\, _4\phi_3\left[\matrix{
-\ep q^{-k}/t,\, -\ep,\, ad,\, \ep^2 uv/q^2\cr
adq^{1-k}/t,\, u\ep^2/q,\, v\ep^2/q\cr} ;q,q\right]\cr
&+{(\ep^2, ad, bctu/q, bctv/q, uv\ep^2/q^2;q)_\ty\o
(bc, ad/t, u\ep^2/q, v\ep^2/q, bcuvt/q^2;q)_\ty}\cr
&\cdot \tk {(t, -\ep/ad, -\ep t/ad, bcuvt/q^2;q)_k\o
(q, qt/ad, bctu/q, bctv/q;q)_k}\,q^k\cr
&\qquad\qquad \cdot\, _4\phi_3\left[\matrix{
q^{-k},\, -t,\, tq^{1/2},\, -tq^{1/2}\cr
qt^2,\, -\ep t/ad,\, -adq^{1-k}/\ep\cr} ;q,q\right].\cr
}
$$
The contribution of the first term in (4.6) to $K_t^{\lm,\mu}
(x,y)$
is
$$\leqalignno{
K_t^{(1)} (x,y) &= B(\ta, \phi)\, (1-t^2)\, {(-qt\ep ;q)_\ty\o
(-t/\ep;q)_\ty}&(4.7)\cr
&\cdot \tk {(\ep, \ep q^{1/2}, -\ep q^{1/2}, -\ep/ad;q)_k\o
(q, bc, -qt\ep, -q\ep/t;q)_k}\, q^k\cr
&\cdot \sum^k_{\ell = 0} {(q^{-k}, -\ep, ad;q)_\ell\o
(q, -adq^{1-k}/\ep ;q)_\ell}\, q^\ell \, U_\ell,\cr
}
$$
where
$$\leqalignno{
U_\ell = &  \int^{q\meit/b}_{q\eit/b} d_q u \,{(bu\eit, bu\meit,
u\ep^2
q^{\ell -1};q)_\ty\o
(bau/q, bcu/q, bdu/q;q)_\ty} &(4.8)\cr
&\cdot \int^{q\meip/\g}_{q\eip/\g}d_q v \, {(\g v\eip, \g v\meip,
v\ep^2
q^{\ell -1}, uv\ep^2/q;q)_\ty\o
(\g\a u/q, \g\b v/q, \g\d v/q, uv\ep^2 q^{\ell-2};q)_\ty}.\cr
}
$$
\indent
Likewise, the contribution of the second term in (4.6) is
$$\leqalignno{
K^{(2)}_t (x,y) &= B(\ta, \phi)\, {(\ep^2, -\ep/ad, t, -
t\ep/ad;q)_\ty\o
(-\ep, bc, t/ad, -\ep/t;q)_\ty}&(4.9)\cr
&\cdot \tk {(-t, tq^{1/2}, -tq^{1/2}, t/ad;q)_k\o
(q, qt^2, -t\ep/ad, -qt/\ep;q)_k} \, q^k\cr
&\cdot \sum^\ty_{\ell = 0} {(-\ep q^{-k}/t, -\ep, ad;q)_\ell\o
(q, adq^{1-k}/t;q)_\ell}\, q^\ell \, U_\ell.\cr
}
$$
Finally, the last term on the right side of (4.6) gives
$$\leqalignno{
K^{(3)}_t (x,y)&= B(\ta, \phi) \, {(\ep^2, ad;q)_\ty\o
(bc, ad/t;q)_\ty}& (4.10)\cr
&\cdot \tk {(t, -\ep/ad, -t\ep/ad;q)_k\o
(q, qt/ad;q)_k} \, q^k\cr
&\cdot\, _4\phi_3\left[\matrix{
q^{-k}, -t, tq^{1/2}, -tq^{1/2}\cr
qt^2, -t\ep/ad, -adq^{1-k}/\ep\cr} ;q,q\right] V_k,\cr
}
$$
where
$$\leqalignno{
V_k = &  \int^{q\meit/b}_{q\eit/b} d_q u \, {(bu\eit, bu\meit,
bctuq^{k-1};q)_\ty\o
(bau/q, bcu/q, bdu/q;q)_\ty} &(4.11)\cr
&\cdot \int^{q\meip/\g}_{q\eip/\g}d_q v \, {(\g v \eip, \g v\meip,
bctvq^{k-1}, uv\ep^2/q^2;q)_\ty\o
(\g\a v/q, \g\b v/q, \g\d v/q, bcuvtq^{k-2};q)_\ty}.\cr
}
$$
{\bf 5. Computation of $\bf K^{(1)}_t (x,y)$ and $\bf
K_t^{(2)}(x,y)$}.\quad Using [10, (2.10.19)], we get
$$\leqalignno{
\int^{q\meip/\g}_{q\eip/\g} & {(\g v \eip, \g v\meip, v\ep^2
q^{\ell-1}, uv\ep^2/q;q)_\ty\o
(\g\a v/q, \g\b v/q, \g\d v/q, uv\ep^2q^{\ell-2};q)_\ty} \,
d_q v&(5.1)\cr
&= {q(1-q)\o 2i\g} \, (q, \a\b, \a\d, \b\d;q)_\ty \, h (y;\g)
\,\rho^\mu
(y)\cr
&\cdot {(\a\eip, \a\meip,\a\b u/q;q)_\ell\o
(\a\b, \a\d, \a/\d;q)_\ty}\cr
&\cdot\ _8W_7 \left(\d q^{-\ell}/\a; \d\eip, \d\meip, \b\d u/q,
q^{1-\ell}/\a\b, q^{-\ell};q, q^2/\a\d u\right).\cr
}
$$
Substituting (5.1) in (4.8) we get
$$\leqalignno{
U_\ell &= {q(1-q)\o 2i\g}\, (q, \a\b, \a\d, \b\d;q)_\ty  \, h
(y;\g)\,
\rho^\mu (y) \, {(\a\eip, \a\meip;q)_\ell\o
(\a\b, \a\d, \a/\d;q)_\ell}&(5.2)\cr
&\cdot\sum^\ell_{m=0} {(1-\d q^{2m-\ell}/\a) (\d q^{-\ell}/\a,
\d\eip, \d\meip, q^{1-\ell}/\a\b, q^{-\ell};q)_m\o
(1-\d q^{-\ell}/\a) (q, q^{1-\ell} \meip/\a, q^{1-\ell}\eip/\a,
\b\d,
q\d/\a;q)_m} \, I_m,\cr
}
$$
where
$$\leqalignno{&&(5.3)\cr
I_m &= q^{\ell m - {m\choose 2}} (-\b/\d)^m
\int^{q\meit/b}_{q\eit/b} {(bu\eit, bu\meit, u\ep^2 q^{\ell-1},
\a\b u/q;q)_\ty\o
(bcu/q, bduq^{m-1}, \a\b uq^{\ell-m-1}, abu/q;q)_\ty} \, d_q u \cr
&= q^{\ell m - {m\choose 2}} (-\b/\d)^m \, {q(1-q)\o b} \, \meit
\cr
&\cdot{(q, e^{2i\ta}, qe^{-2i\ta}, \a\b dq^\ell/b, \a\b cq^{\ell-
m}/b, cdq^m, q^\ell
\ep^2\meit/b, \a\b \meit/b;q)_\ty\o
(c\eit, dq^m\eit, \a\b q^{\ell-m}\eit/b, c\meit, dq^m\meit, \a\b
q^{\ell-m}\meit/b, a\meit, \a\b cdq^\ell \meit/b;q)_\ty}\cr
&\cdot\, _8W_7\left(\a\b cdq^{\ell-1}\meit/b; c\meit, dq^m\meit,
\a\b
q^{\ell-m}\meit/b, cdq^\ell, \a\b/ab;q, a\eit \right),\cr
}
$$
by [10, (2.10.19)]. Simplifying the coefficients, we obtain from
(5.2), (5.3) and (4.5)
$$\leqalignno{
U_\ell &= {B^{-1} (\ta, \phi) \, (\a\b c/b, \a\b d/b, a\eit,
\ep^2\meit/b;q)_\ty\o
(ac, ad, \a\b \eit/b, \a\b cd\meit/b;q)_\ty}&(5.4)\cr
&\cdot {(\a\eip, \a\meip, \a\b\eit/b, \a\b\meit/b, \a\b
cd\meit/b;q)_\ell\o
(\a\b, \a\d, \a/\d, \a\b d/b, \a\b c/b, \ep^2\meit/b;q)_\ell}\cr
&\cdot \sum^\ell_{m=0} {(1-\d q^{2m-\ell}/\a) (\d q^{-\ell}/\a,
q^{1-\ell}/\a\b, bq^{1-\ell}/\a\b c, q^{-\ell};q)_m\o
(1-\d q^{-\ell}/\a) (q, \b\d, cd, q\d/\a;q)_m}\cr
&\cdot {(\d\eip, \d\meip, d\eit, d\meit;q)_m\o
(q^{1-\ell}\meip/\a, q^{1-\ell}\eip/\a, bq^{1-\ell}\meit/\a\b,
bq^{1-\ell}\eit/\a\b;q)_m} (bcq/\a\d)^m\cr
&\cdot\,_8W_7\left(\a\b cdq^{\ell-1}\meit/b;c\meit, dq^m\meit, \a\b
q^{\ell -m}\meit/b, cdq^\ell, \a\b/ab;q, a\eit\right).\cr
}
$$
Use of (5.4) in (4.7) and (4.9) then gives (3.10) and (3.11).
\bigskip
\noindent
{\bf 6. Computation of $\bf K^{(3)}_t (x,y)$}.\quad In order to
compute $K^{(3)}_t (x,y)$, we have to express $V_k$ of (4.11) in a
form that is real in both $\ta$ and $\phi$. First of all, by [10,
(2.10.19)]
$$\leqalignno{
I(u) &\colon = \int^{q\meip/\g}_{q\eip/\g} {(\g v\eip, \g v\meip,
bct vq^{k-1}, uv\ep^2/q;q)_\ty\o
(\g\a v/q, \g\b v/q, \g\d v/q, bcuvtq^{k-2};q)_\ty} \, d_qv
&(6.1)\cr
&= {q(1-q)\o 2i\g} \, (q, \a\b, \a\d, \b\d;q)_\ty \, h(y;\g) \,
\rho^\mu
(y)\cr
&\cdot {(bctq^k\meip/\g, u\ep^2\meip/\g q;q)_\ty\o
(bcutq^{k-1}\meip/\g, \ep^2\meip/\g;q)_\ty}\cr
&\cdot\, _8W_7\bg{ \ep^2\meip/\g q; \a\meip, \b\meip, \d\meip, q/u,
adq^{-k}/t;q, bcutq^{k-1}\eip/\g}.\cr
}
$$
Using Bailey's transformation formula [10, (2.10.1)] for a
very-well-poised $_8\phi_7$ series we find that
$$\leqalignno{
I(u) &= {q(1-q)\o 2i\g} \, (q, \a\d, \b\d;q)_\ty \, h(y;\g) \,
\rho^\mu
(y)&(6.2)\cr
&\cdot {(bctq^k\eip/\g, bctq^k\meip/\g, \a\b u/q, bc\d utq^{k-
1}/\g;q)_\ty\o
(bcutq^{k-1}\eip/\g, bcutq^{k-1} \meip/\g, bc\d tq^k/\g;q)_\ty}\cr
&\cdot\, _8W_7\bg{ bc\d tq^{k-1}/\g; \d\eip, \d\meip, \a\d tq^k/ad,
bctq^k/\a\g, q/u;q, \a\b u/q}.\cr
}
$$
Now we break up the $_8W_7$ series in (6.2) into two balanced
nonterminating $_4\phi_3$ series by use of [10, (2.10.10)]:
$$\leqalignno{
I(u)&= {(\a\eip, \a\meip, \b\d, \a\b u/q, bctq^k/\g\d, bc\d utq^{k-
1}/\g;q)_\ty\o
(\a/\d, bcutq^{k-1}\eip/\g, bcutq^{k-1}\meip/\g;q)_\ty} &(6.3)\cr
&\cdot\, _4\phi_3\left[\matrix{
bctq^k/\a\g,\, \d\eip,\, \d\meip, \,\b\d u/q \cr
q\d/\a,\, \b\d, \,bc\d utq^{k-1}/\g\cr};q,q\right] + \ {\rm idem}
(\a;\d).\cr
}
$$
Use of (6.3) in (4.11) now gives
$$\leqalignno{
V_k &= {q(1-q)\o 2i\g} \, (q;q)_\ty  \, h(y;\g) \, \rho^\mu (y) \,
{(\a\eip,
\a\meip, \b\d, bctq^k/\g\d;q)_\ty\o (\a/\d;q)_\ty}&(6.4)\cr
&\cdot \sum^\ty_{\ell = 0} {(bctq^k/\a\g, \d\eip, \d\meip;q)_\ell\o
(q, \b\d, q\d/\a;q)_\ell} \, q^\ell\cr
&\cdot \int^{q\meit/b}_{q\eit/b} {(bu\eit, bu\meit, bctuq^{k-1},
bc\d tuq^{k+\ell-1}/\g;q)_\ty\o
(bcu/q, bcutq^{k-1}\eip/\g, bcutq^{k-1}\meip/\g,
\b\d uq^{\ell-1};q)_\ty}\cr
&\cdot {(\a\b u/q, \d\b u/q;q)_\ty\o
(abu/q, bdu/q;q)_\ty} \, d_q u\cr
&+\ {\rm idem} (\a;\d).\cr
}
$$
Observe that the two $q$-integrals on the right side of (6.4)
differ only in that $\a$ and $\d$ are interchanged. So we need to
transform only one of them. By using the definition (4.3) we obtain
$$\leqalignno{
\int^{q\meit/b}_{q\eit/b} &{(bu\eit, bu\meit, bctuq^{k-1}, bc\d
tuq^{k+\ell-1}/\g;q)_\ty\o
(bcu/q, bcutq^{k-1}\eip/\g, bcutq^{k-1}\meip/\g, \b\d uq^{\ell-
1};q)_\ty}&(6.5)\cr
&\cdot {(\a \b u/q, \d\b u/q;q)_\ty\o (abu/q, bdu/q;q)_\ty} \, d_q
u\cr
&= {q(1-q)\o 2ib} \, (q;q)_\ty \, h(x;b) \, \rho^\lm (x) \,
{(\b\d\eit/b;q)_\ell\o (c\d t\eit q^k/\g;q)_\ell}\cr
&\cdot {(ctq^k\eit, c\d tq^k\eit/\g, a\meit, c\meit, d\meit,
\a\b\eit/b;q)_\ty\o
(ctq^k e^{i(\ta + \phi)}/\g, ctq^k e^{i(\ta-\phi)}/\g, e^{-
2i\ta};q)_\ty}\cr
&\cdot\, _6\phi_5\left[\matrix{
ctq^k e^{i(\ta + \phi)}/\g, ctq^k e^{i(\ta -\phi)}/\g, \d\b q^\ell
\eit/b, a\eit, c\eit, d\eit\cr
qe^{2i\ta}, ctq^k\eit, c\d tq^{k+\ell}\eit/\g, \a\b \eit/b,
\d\b\eit/b\cr} ;q,q\right]\cr
&- \ {\rm idem}\ (\ta; -\ta).\cr
}
$$
Using (6.5) and (4.5) in (6.4) we find that
$$\leqalignno{
V_k &= B^{-1} (\ta, \phi) \, \Biggl\{  {(\a\eip, \a\meip, bct/\a\d,
ct\eit, c\d t\eit/\g;q)_\ty\o
(ac, ad, cd, \a\b, \a\d, \a/\d;q)_\ty}&(6.6)\cr
&\cdot{(a\meit, c\meit, d\meit, \a\b\eit/b;q)_\ty (cte^{i(\ta +
\phi)}/\g, cte^{i(\ta-\phi)}/\g;q)_k\o
(cte^{i(\ta + \phi)}/\g, cte^{i(\ta -\phi)}/\g;q)_\ty (bct/\g\d,
ct\eit, c\d t\eit/\g;q)_k}\cr
&\cdot \sum^\ty_{\ell = 0} {(bctq^k/\a\g, \d\eip, \d\meip,
\b\d\eit/b;q)_\ell\o
(q, q\d/\a, \b\d, c\d tq^k\eit/\g;q)_\ell} \, q^\ell\cr
&\cdot\, _6\phi_5\left[\matrix{
ctq^ke^{i(\ta + \phi)}/\g, ctq^k e^{i(\ta-\phi)}/\g, \b\d
q^\ell\eit/b, a\eit, c\eit, d\eit\cr
qe^{2i\ta}, ctq^k\eit, c\d tq^{k+\ell}\eit/\g, \a\b\eit/b,
\d\b\eit/b\cr};q,q\right]\cr
&+\ {\rm idem} (\ta; -\ta)\Biggr\} + \ {\rm idem} (\a;\d).\cr
}
$$
Changing the order of summation one can  rewrite the sum over
$\ell$
in (6.6) in the following form
$$
\eqalign
{
&\cdot \sum^\ty_{\ell = 0} {(ctq^ke^{i(\ta + \phi)}/\g,
ctq^ke^{i(\ta - \phi)}/\g, a\eit, c\eit, d\eit;q)_\ell \o
(q, qe^{2 i\ta }, ctq^k\eit, \a\b\eit/b, c\d tq^k\eit/\g;q)_\ell}
\, q^\ell\cr
&\quad \cdot\, _4\phi_3\left[\matrix{
bctq^k/\a\g, \, \d\eip, \, \d\meip, \, \d \b q^\ell \eit /b \cr
q\d/\a, \, \b\d, \, c\d tq^{k+\ell}\eit/\g \cr}; \,q, \,q
\right],\cr
}
$$
where the $\,_4\phi_3$ function is now balanced. The two 
$\,_4\phi_3$
series that arise from (6.6) in this way can be combined into a
single
$\,_8\phi_7$ series by [10, III.36].
Substituting this result in (4.10) we get (3.12).
\bigskip
\noindent
{\bf 7. Multiplication law for the kernels}.\quad It follows from
(2.4), (3.7) and (3.8) that
$$\int^1_{-1} P_t^{\lm,\mu} (x,y) \, P_{t'}^{\mu,\lm'} (y, x') \,
\rho^\mu (y)\, dy = P_{tt'}^{\lm,\lm'} (x, x'),\leqno (7.1)$$
when $\max (|\lm|, |\mu|, |\lm'|, |t|, |t'|)< 1$, where $\lm = (a,
b, c, d)$, $\lm' = (a', b', c', d')$, $\mu = (\a, \b, \g, \d)$, and
$|\lm| < 1$ means $|a|, |b|, |c|, |d|$ are all numerically less
than 1. The corresponding formula for the kernels $K_t^{\lm,\mu}
(x,y)$ is
$$\leqalignno{
&\int^1_{-1} K_t^{\lm, \mu} (x,y) \, K_{t'}^{\mu, \lm'} (y, x')\,
\rho^\mu (y) \, dy&(7.2)\cr
&\quad= {2\pi (\a\b\g\d;q)_\ty\o
(q, \a\b, \a\g, \a\d, \b\g, \b\d, \g\d;q)_\ty}\,
K_{tt'}^{\lm, \lm'}(x, x').\cr
}
$$
\indent
In the next sections we will consider an analytic continuation
of (7.2) when the parameters $t$ and $t'$ go to the unit circle and
$tt'$ approaches 1, in order to obtain the orthogonality relation
for the kernels $K_t^{\lm,\mu}(x,y)$.
\bigskip
\noindent
{\bf 8. Computation of ${\bf K_1^{\lm, \mu} (x,y)}$}.\quad
Setting $t = 1$ in (3.10) and (3.11) we can see that $K^{(1)}_1
(x,y) = 0 = K_1^{(2)} (x,y)$. So, from (3.9), (3.12) and (4.10) we
find that
$$\leqalignno{
K_1^{\lm,\mu} (x,y) &= K^{(3)}_1 (x,y)&(8.1)\cr
&= B(\ta, \phi) \, {(\ep^2;q)_\ty\o
(bc;q)_\ty} \, V_0\Big|_{t=1}.\cr
}
$$
However, from (6.4) we have
$$\leqalignno{
V_0\Big|_{t=1} &= {q(1-q)\o 2i\g} \, (q;q)_\ty \, h(y;\g) \,
\rho^\mu
(y)&(8.2)\cr
&\cdot \Biggl[ {(\a\eip, \a\meip, \b\d, bc/\g\d;q)_\ty\o
(\a/\d;q)_\ty} \, \sum^\ty_{\ell = 0} {(bc/\a\g, \d\eip,
\d\meip;q)_\ell\o
(q, q\d/\a, \b\d;q)_\ell} \, q^\ell\cr
&\cdot \int^{q\meit/b}_{q\eit/b} {(bu\eit, bu\meit, bc\d uq^{\ell-
1}/\g, \a\b u/q;q)_\ty\o
(bcu\eip/\g q, bcu\meip/\g q, \b\d uq^{\ell-1}, abu/q;q)_\ty} \,
d_q u\cr
&+ {(\d\eip, \d\meip, \a\b, bc/\a\g;q)_\ty \o
(\d/\a;q)_\ty} \, \sum^\ty_{\ell = 0} {(bc/\g\d, \a\eip,
\a\meip;q)_\ell\o
(q, q\a/\d, \a\b;q)_\ell} \, q^\ell\cr
&\cdot \int^{q\meit/b}_{q\eit/b} {(bu\eit, bu\meit, bc\a uq^{\ell
-
1}/\g, \a\b u/q;q)_\ty\o
(bcu\eip/\g q, bcu\meip/\g q, \a\b uq^{\ell-1}, abu/q;q)_\ty} \,
d_q u
\Biggr].\cr
}
$$
\indent
All of our calculations so far have been done on the basis of only
two conditions (3.3) connecting the $\lm$-and $\mu$-parameters. We
shall now impose a third condition:
$$\b\g = bc.\leqno (8.3)$$
This still leaves us with enough freedom to compute $K_1^{\lm, \mu}
(x,y)$ without hitting any singularities while enabling us to
simplify the computations in (8.2) enormously. First observe that
the $\ell$-dependent terms in both $q$-integrals in (8.2) cancel
out so that right away one has a simpler formula
$$\leqalignno{
V_0\Big|_{t=1} &= {q(1-q)\o 2i\g} \, (q;q)_\ty \, h(y;\g) \,
\rho^\mu
(y)&(8.4)\cr
&\cdot \int^{q\meit/b}_{q\eit/b} {(bu\eit, bu\meit, \a\b
u/q;q)_\ty\o
(\b u\eip/q, \b u\meip/q, abu/q;q)_\ty} \, d_q u\cr
&\cdot \Biggl\{ {(\a\eip, \a\meip, \b\d, \b/\d;q)_\ty\o
(\a/\d;q)_\ty} \ _3\phi_2\left[\matrix{
\b/\a, \d\eip, \d\meip\cr
q\d/\a, \b\d\cr} ;q,q\right]\cr
&+ {(\d\eip, \d\meip, \a\b, \b/\a;q)_\ty\o
(\d/\a;q)_\ty}\ _3\phi_2\left[\matrix{
\b/\d, \a\eip, \a\meip\cr
q\a/\d, \a\b\cr} ;q;q\right] \Biggr\}.\cr
}
$$
By the nonterminating $_3\phi_2$ summation formula [10, II.24] we
can sum the expression inside the curly brackets above to get
$$\leqalignno{
V_0\Big|_{t =1} &= {q(1-q)\o 2i\g} \, (q;q)_\ty \, h(y;\g) \,
\rho^\mu (y)\,
(\a\d, \b\eip, \b\meip;q)_\ty&(8.5)\cr
&\cdot \int^{q\meit/b}_{q\eit/b} {(bu\eit, bu\meit, \a\b
u/q;q)_\ty\o
(\b u\eip/q, \b u\meip/\g, abu/q;q)_\ty} \, d_q u.\cr
}
$$
Because of (3.3) and (8.3), the $q$-integral in (8.5) can also be
evaluated by [10, (2.10.18)]. Our final result for the kernel in
(8.1) is
$$\leqalignno{
K_1^{\lm, \mu} (x,y) &= {(abcd;q)_\ty\o
(\a\b, ac, ad, bc, bd, cd;q)_\ty}&(8.6)\cr
&\cdot {(\a\eip, \a\meip, \b\eip, \b\meip, c\eit, c\meit, d\eit,
d\meit, (\b/b)^2;q)_\ty\o
(\b e^{i\ta + i\phi}/b, \b e^{i\ta -i\phi}/b, \b e^{i\phi- i\ta}/b,
\b e^{-i\ta - i\phi}/b;q)_\ty},\cr
}
$$
where, it has been assumed that $|\b| < |b|$.
\medskip
We would like to point out that it is possible to derive (8.6)
directly
from (3.12) by using [10, II.24].
\bigskip
\noindent
{\bf 9. Relation of ${\bf K_1^{\lm, \mu} (x,y)}$ with the
continuous $q$-Hermite polynomials}.\quad The Poisson kernel for
the
continuous $q$-Hermite polynomials is, see [3],
$$\leqalignno{
&\sn H_n(\cos\ta |q) \, H_n (\cos\phi|q) \, {r^n\o (q;q)_n}
&(9.1)\cr
&\quad = {(r^2;q)_\ty\o
(re^{i(\ta + \phi)}, re^{i(\ta-\phi)}, re^{i(\phi-\ta)}, re^{-i(\ta
+ \phi)};q)_\ty}, \quad |r| < 1.\cr
}
$$
So the right side of (8.6) coincides with it (up to a factor
independent of $r$) if we identify $r$ with $\b/b$. The $q$-Hermite
polynomials are orthogonal on the compact interval $[-1, 1]$ and so
are complete in $L^2 [-1, 1]$, [19, Thm. 3.1.5]. We also have
$$\lim_{r\to 1^-} \int^1_{-1} K^{(0)}_r (x,y) f(y) \, dy = f(x),
\quad
a.e.,\leqno (9.2)$$
for every function $f$ continuous on $[-1, 1]$, where
$$K^{(0)}_r (x,y) = {(q, r^2, e^{2i\phi}, e^{-2i\phi};q)_\ty\o
2\pi \sin\phi (re^{i(\ta + \phi)}, re^{i(\ta - \phi)}, re^{i(\phi
-\ta)}, re^{-i(\ta + \phi)};q)_\ty}.\leqno (9.3)$$
It also means that in the limit $r \to 1^-$ the kernel
$K^{(0)}_r(x,y)$ behaves like the delta functional $\d (x-y)$.
\medskip
The proof of (9.2) follows by Wiener's arguments in lemma [21,
proposition $X_{55}$, p.~41]. Indeed, if we set $a = 0,\ b = \pi$
and
$K_r(\ta, \phi) = \sin\phi \, K^{(0)}_r (\cos\ta, \cos\phi)$,
then the hypotheses of this lemma are satisfied by this kernel.
Here
$$\leqalignno{
\int^\pi_0 K_r (\ta, \phi) \, d\phi &= \sn {r^n\o (q;q)_n} \,
H_n(\cos\ta |q)&(9.4)\cr
&\cdot {(q;q)_\ty\o 2\pi} \int^\pi_0 H_n (\cos\phi| q) \,
(e^{2i\phi}, e^{-2i\phi}; q)_\ty \, d\phi= 1\cr
}
$$
due to uniform convergence of the series for $|r| < 1$ and the
orthogonality of the $q$-Hermite polynomials. Thus the property
(4.12) of [21] is established. Again, if $\ta + \phi \ge |\ta -
\phi| > \psi > 0$, then
$$\leqalignno{
K_r (\ta, \phi) &< {(q, r^2, e^{2i\phi}, e^{-2i\phi};q)_\ty\o
2\pi (r e^{i\psi}, re^{-i\psi};q)^2_\ty} &(9.5)\cr
&<\ {\rm (constant)}\ {(r^2;q)_\ty\o
(re^{i\psi}, re^{-i\psi};q)_\ty^2},\cr
}
$$
since
$$\eqalignno{
&\bg{re^{i(\ta\mp \phi)}, re^{-i(\ta\mp \phi)};q}_\ty\cr
& = \prod^\ty_{k=0} (1- 2r \cos (\ta\mp\phi) q^k + r^2 q^{2k})\cr
& > \prod^\ty_{k=0} (1 - 2r \cos \psi q^k + r^2 q^{2k}) =
(re^{i\psi}, re^{-i\psi};q)_\ty.\cr
}
$$
The validity of [21, (4.13)] is thus assured. The positivity of
$K_r (\ta, \phi)$ is also obvious. Therefore, by [21, proposition
$X_{55}$]
$$\lim_{r\to 1^-} \int^\pi_0 K_r (\ta, \phi) \, g(\phi) \, d\phi =
g(\ta), \quad a.e.,\leqno (9.6)
$$
which is equivalent to (9.2).
\bigskip
\noindent
{\bf 10. Complex orthogonality of the kernels}.\quad Set $t =
re^{i\tau}, t' = re^{-i\tau}$ with $0 < \tau < \pi$ in eq. (7.2)
and consider its limiting form when $r \to 1^-$ and then $\lm' \to
\lm$.  With the aid of (8.6), (9.2) and (9.3) we can write (in the
distribution sense), for the right side of (7.2):
$$\leqalignno{
K_1^{\lm,\lm'} (x, x') \to \  & {2\pi (abcd;q)_\ty\o
(q, ab, ac, ad, bc, bd, cd; q)_\ty}&(10.1)\cr
&\cdot {\d (x-x')\o \rho^\lm (x)},\quad {\rm as}\ \lm'\to\lm.\cr
}
$$
Here $\lm = (a,b,c,d)$, $\lm' = (a', b', c', d')$ with $ac = a'c'$,
$bd = b'd'$, $bc = b'c'$ and $b' < b$.
\medskip
To consider the analytic continuation of the left side of (7.2)
when $t = (t')^* = re^{i\tau} = e^{i\tau -\kappa}$ and $\kappa \to
0^+$,
let us
rewrite this formula in terms of a contour integral (Fig.~1):
$$\leqalignno{
&\int^1_{-1} K^{\lm, \mu}_{re^{i\tau}} (x,y) \,
K^{\mu, \lm'}_{re^{-i\tau}} (y, x') \, \rho^\mu (y) \, dy&(10.2)\cr
&\quad = {\log q^{-1}\o i} \int_{C^-} {L_{re^{i\tau}}^{\lm, \mu}
(z, s) \,
L_{re^{-i\tau}}^{\mu, \lm'} (s, z') \, (q^{2s}, q^{-2s};q)_\ty \,
ds \o
(\a q^s, \a q^{-s}, \b q^s, \b q^{-s}, \g q^s, \g q^{-s}, \d q^s,
\d q^{-s};q)_\ty},\cr
}
$$
where $x = {1\o 2} (q^z + q^{-z}), q^z = \eit$; $y = {1\o 2} (q^s
+ q^{-s})$, $q^s =\eip$, and $L_t^{\lm, \mu} (z, s) = K_t^{\lm,
\mu} (x,y)$. An examination of the series on the right sides of
(3.10)--(3.12) reveals that the poles only originate from the
function $K^{(3)}_t (x,y)$ given by (3.12). Even here the only
singular term corresponds to $k = 0$ on the right side of (3.12).
The poles of the integrand which are located in the lower half
plane,
are given by
$$\eqalign{
s_1 &= {\om + \kappa - i (\tau + \ta)\o
\log q^{-1}},\quad s_1 + 1, \ldots,\cr
s_2 &= {\om + \kappa - i(\tau -\ta)\o
\log q^{-1}},\quad s_2 + 1, \ldots\ ;\cr} \leqno(10.3)
$$
and
$$\eqalign{
s_1^\prime &= {\om' - \kappa - i (\tau + \ta')\o
\log q^{-1}},\quad s_1^\prime - 1, \ldots,\cr
s_2^\prime &= {\om' - \kappa - i(\tau -\ta')\o
\log q^{-1}},\quad s_2^\prime - 1, \ldots\ .\cr} \leqno(10.4)
$$
Here $c/\g = e^{-\om}, c'/\g = e^{-\om'}, r = e^{-\kappa}$ and $c
< c' < \g$. In the limit $\kappa \to 0^+$ the first poles of
(10.4), $s^\prime_1$ and $s^\prime_2$, will move from the left half
plane to the right plane,  thus invalidating the formula (10.2).
Therefore, we have to replace our original contour $C^-$ by a
contour $C^-_\mu$ with certain indentations (see Fig.~2), thus
separating the increasing and decreasing sequences of poles (10.3)
and (10.4). After the analytic continuation of both sides of (7.2)
we arrive at the complex orthogonality property of the kernels:
$$\leqalignno{
{\log q^{-1}\o i}
& \int_{C^-_\mu} {L^{\lm, \mu}_{e^{i\tau}} (z, s) \,
L^{\mu,\lm}_{e^{-i\tau}} (s, z';) (q^{2s}, q^{-2s};q)_\ty \, ds\o
(\a q^s, \a q^{-s}, \b q^s, \b q^{-s}, \g q^s, \g q^{-s}, \d q^s,
\d q^{-s};q)_\ty}&(10.5)\cr
&= {2\pi (\a\b\g\d;q)_\ty\o
(q, \a\b, \a\g, \a\d, \b\g, \b\d, \g\d;q)_\ty}\cr
&\cdot {2\pi (abcd;q)_\ty\o
(q, ab, ac, ad, bc, bd, cd; q)_\ty} \cdot {\d (x-x')\o
\rho (x; a, b, c, d)}.\cr
}
$$
For the $P$-kernels, the orthogonality relation (10.5) has the more
compact form
$$\int_{\G_\mu} P^{\lm, \mu}_{e^{i\tau}} (x,y) \, P^{\mu,
\lm}_{e^{-
i\tau}} (y, x') \, \rho^{\mu} (y) \, dy = {\d (x-x')\o \rho^\lm
(x)}.\leqno (10.6)$$
Here $C^-_\mu \to \G_\mu$ when $y = {1\o 2} (q^s + q^{-s})$.
%\bigskip
\vfill \eject
\noindent
{\bf 11. Real orthogonality of the kernels: the case of one
additional mass-point}.\ The structure of poles in (10.3) and
(10.4) of the integrand in (10.5) was shown in Fig.~2, where we now
let $\om' \to \om$ (or $c' \to c)$. First we consider the simplest
case when $0 < \om/\log q^{-1} < 1$ or $1 < \g/c < q^{-1}$, so only
the first pole from each sequence has appeared in the lower half
$s$-plane after analytic continuation. Let us denote the integrand
on the left side of (10.5) as $F(z, z', s)$. Then
$$\leqalignno{
\int_{C^-_\mu} F(z, z', s) \, ds &= \int_{C^-} F(z, z', s) \,
ds&(11.1)\cr
&+ 2\pi i \left(  Res\ F(z, z', s)\Big|_{s=s^\prime_1} + \ Res\
F(z,
z', s)\Big|_{s= s^\prime_2}  \right), \cr
}
$$
where
$$Res\ F(z, z', s)\Big|_{s = s^\prime_\a} = \lim_{s\to s^\prime_\a}
(s-s^\prime_\a) \, F(z, z', s), \ \a = 1, 2.\leqno (11.2)$$
In view of (3.9)--(3.12) we get
$$\leqalignno{
Res\ F(z, z', s)& \Big|_{s = s^\prime_1}&(11.3)\cr
&= {L^{\lm, \mu}_{e^{i\tau}} (z, s) (q^{2s}, q^{-2s};q)_\ty\o
(\a q^s, \a q^{-s}, \b q^s, \b q^{-s}, \g q^s, \g q^{-s}, \d q^s,
\d q^{-s};q)_\ty} \Bigg|_{s = s^\prime_1}\cr
&\cdot \lim_{s\to s^\prime_1} (s- s^\prime_1) \,
L^{\mu,\lm}_{e^{-i\tau}} (s, z'),\cr
}
$$
and
$$\leqalignno{
&\lim_{s\to s^\prime_1} (s-s^\prime_1) \,
L^{\mu,\lm}_{e^{-i\tau}} (s, z')&(11.4)\cr
&= {(\ep^2, aq^{z'}, cq^{z'}, dq^{z'}, \a q^{-s^\prime_1}, \g q^{-
s^\prime_1}, \d q^{-s^\prime_1};q)_\ty\o
\log q^{-1} (q, \a\d e^{i\tau}, \a\g, \b\g, \d\g, ad, q^{-
2s^\prime_1}, q^{2z'};q)_\ty}\cr
&\cdot \Biggl\{ {(\b\g e^{-i\tau}/cd, aq^{-z'},
abq^{s^\prime_1}/\b;q)_\ty\o
(ab, a/d;q)_\ty} \ _3\phi_2\left[\matrix{
\b\g e^{-i\tau}/ac, dq^{-z'}, bdq^{s^\prime_1}/\b\cr
qd/a, bd\cr} ;q,q\right]\cr
&+ {(\b\g e^{-i\tau}/ac, dq^{-z'}, bdq^{s^\prime_1}/\b;q)_\ty\o
(bd, d/a;q)_\ty}\ _3\phi_2\left[\matrix{
\b\g e^{-i\tau}/cd, aq^{-z'}, abq^{s^\prime_1}/\b\cr
qa/d, ab\cr} ;q,q\right] \Biggr\},\cr
}
$$
where $q^{s^\prime_1} = cq^{z'} e^{i\tau}/\g$. We can now sum the
expression above by [10, II.24], thus getting the final result
\vfill
\eject
$$\leqalignno{&&(11.5)\cr
&\lim_{s \to s^\prime_1} (s-s^\prime_1) \,
L^{\lm,\mu}_{e^{-i\tau}}(s, z')\cr
&= {(\ep^2, \a\g e^{-i(\tau + \ta')}/c, \b\g e^{-i(\tau + \ta')}/c,
\g^2 e^{-i(\tau + \ta')}/c, \d\g e^{-i(\tau + \ta')}/c, ae^{i\ta'},
be^{i\ta'}, ce^{i\ta'}, de^{i\ta'};q)_\ty\o
\log q^{-1} (q, \a\g, \b\g, \d\g, ab, ad, bd, e^{2i\ta'}, \g^2 e^{-
2i(\tau + \ta')}/c^2;q)_\ty},\cr
}
$$
since $q^{z'} = e^{i\ta'}$ and $q^{s^\prime_1} = ce^{i\tau}
q^{z'}/\g$. The value of $Res\ F(z, z', s)\Big|_{s = s^\prime_2}$
is obtained from (11.3) and (11.5) by simply replacing $\ta'$ by $-
\ta'$. Thus we have shown that
$$\leqalignno{
&\int^1_{-1} P^{\lm,\mu}_{e^{i\tau}} (x,y) \,
P^{\mu,\lm}_{e^{-i\tau}}
(y, x') \rho^\mu (y) \, dy&(11.6)\cr
&+ {(\a\b, \a\d, \b\d;q)_\ty\o
(ab, ad, bd;q)_\ty} \, P^{\lm,\mu}_{e^{i\tau}} (\cos\ta, \cos (\tau
+
\ta' + i\om))\cr
&\cdot {(a e^{i\ta'}, be^{i\ta'}, ce^{i\ta'}, de^{i\ta'},
e^{2i(\tau +\ta') - 2\om};q)_\ty\o
(\a e^{i(\tau + \ta') - \om}, \b e^{i(\tau + \ta') - \om}, \g
e^{i(\tau + \ta') - \om}, \d e^{i(\tau + \ta') - \om},
e^{2i\ta'};q)_\ty}\cr
&+ {(\a\b, \a\d, \b\d;q)_\ty\o
(ab, ad, bd;q)_\ty} \, P^{\lm,\mu}_{e^{i\tau}} (\cos\ta, \cos (\tau
-
\ta' + i\om))\cr
&\cdot {(a e^{-i\ta'}, be^{-i\ta'}, ce^{-i\ta'}, de^{-i\ta'},
e^{2i(\tau - \ta') - 2\om};q)_\ty\o
(\a e^{i(\tau -\ta') - \om}, \b e^{i(\tau -\ta') - \om}, \g
e^{i(\tau
-\ta') - \om}, \d e^{i(\tau -\ta') - \om}, e^{-2i\ta'};q)_\ty}\cr
&= {\d (x - x')\o \rho^\lm (x)}.\cr
}
$$
Here $x = {1\o 2} (e^{i\ta} + e^{-i\ta}) = \cos \ta$, $\g/c =
e^\om$, with $1 < \g/c < q^{-1}$.
\bigskip
\noindent
{\bf 12. Real orthogonality of the kernels: the general case}.\quad
Let us now consider the case of more than one pole, say $N$ poles,
appearing in the lower half $s$-plane which happens when $N -1 <
\om/\log q^{-1} < N$, i.e. $q^{1-N} < \g/c < q^{-N}$, since $e^\om
= \g/c$. Here
$$\leqalignno{
\int_{C^-_\mu} F(z, z', s) \, ds &= \int_{C^-} F(z, z', s)\,
ds&(12.1)\cr
&+ 2\pi i \sum^{N-1}_{j = 0} \left( Res\ F(z, z', s)
\Big|_{s=s^\prime_1} + \ Res\ F(z, z', s)\Big|_{s = s^\prime_2}
\right), \cr
}
$$
where eqs. (11.2) and (11.3) are still valid but with
$q^{s_1^\prime} = {c \o \g} e^{i\tau} q^{z' - j}$ and
$q^{s^\prime_2}
= {c \o \g} e^{i\tau} q^{-z' - j}$, $j = 0, 1, \ldots, N-1$. So,
$$\leqalignno{
&\lim_{s \to s^\prime_1} (s-s^\prime_1) \,
L^{\mu, \lm}_{e^{-i\tau}}(s, z')&(12.2)\cr
&= {(\ep^2, cq^{z'-j}, \a\g e^{-i\tau} q^{j-z'}/c, \g^2 e^{-i\tau}
q^{j-z'}/c, \g\d e^{-i\tau} q^{j-z'}/c;q)_\ty\o
\log q^{-1} (q^{-j};q)_j (q, \a\d e^{i\tau}, \a\g, \b\g, \d\g, ad,
\g^2 e^{-2i\tau} q^{2j - 2z'}/c^2, q^{2z'-j};q)_\ty}\cr
&\cdot \Biggl\{ {(\b\g e^{-i\tau}/cd, aq^{z'}, aq^{-z'}, dq^{z'-j},
abce^{i\tau} q^{z'-j}/\b\g;q)_\ty\o
(ab, a/d;q)_\ty}\cr
&\cdot\, _4\phi_3\left[\matrix{
\b\g e^{-i\tau}/ac, dq^{z'}, dq^{-z'}, bcde^{i\tau} q^{z'-
j}/\b\g\cr
qd/a, bd, dq^{z'-j}\cr} ;q,q\right]\cr
&+ {(\b\g e^{-i\tau}/ac, dq^{z'}, dq^{-z'}, aq^{z'-j}, bcde^{i\tau}
q^{z'-j}/\b\g;q)_\ty\o
(bd, d/a;q)_\ty}\cr
&\cdot\, _4\phi_3\left[\matrix{
\b\g e^{-i\tau}/cd, aq^{z'}, aq^{-z'}, abce^{i\tau} q^{z'-
j}/\b\g\cr
qa/d, ab, aq^{z'-j}\cr} ;q;q\right] \Biggr\}.\cr
}
$$
Both $_4\phi_3$ series above are balanced and their coefficients
are so matched that the sum in the curly brackets can be
transformed to a terminating $_8\phi_7$ series by [10, (2.10.10)]
which, in turn, can be transformed back to a terminating and
balanced $_4\phi_3$. The final result is
$$\leqalignno{&&(12.3)\cr
&\lim_{s\to s^\prime_1}  (s - s^\prime_1) \,
L^{\mu,\lm}_{e^{-i\tau}}(s, z')\cr
&= {(\ep^2, \a\d e^{i\tau} q^{-j}, \a\g e^{-i\tau} q^{j-z'}/c, \b\g
e^{-i\tau} q^{-z'}/c, \g^2 e^{-i\tau} q^{j-z'}/c, \d\g e^{-i\tau}
q^{j-z'}/c;q)_\ty\o
\log q^{-1} (q^{-j};q)_j (q, \a\d e^{i\tau}, \a\g, \b\g, \d\g, ab,
ad, bd;q)_\ty}\cr
&\cdot {(aq^{z'}, bq^{z'}, cq^{z'-j}, dq^{z'};q)_\ty\o
(\g^2 e^{-2i\tau} q^{2j-2z'}/c^2, q^{2z'};q)_\ty}\cr
&\cdot \, _4\phi_3\left[\matrix{
q^{-j}, aq^{-z'}, bq^{-z'}, dq^{-z'}\cr
q^{1-2z'}, \b\g e^{-i\tau} q^{-z'}/c, \a\d e^{i\tau} q^{-j}\cr}
;q,q\right].\cr
}
$$
So, from (3.7), (10.5), (12.1) and (12.3) we obtain the following
general orthogonality property of the kernels:
%
%\vfil\eject
$$\leqalignno{
&\int^1_{-1} P^{\lm,\mu}_{e^{i\tau}} (x,y) \,
P^{\mu,\lm}_{e^{-i\tau}}(y, x') \, \rho^\mu (y) \, dy&(12.4)\cr
&+ {(\a\b, \a\d, \b\d, ae^{i\ta'}, be^{i\ta'}, ce^{i\ta'},
de^{i\ta'}, e^{2i(\ta' + \tau) - 2\om};q)_\ty\o
(ab, ad, bd, \a e^{i(\tau + \ta') -\om}, \b e^{i(\tau + \ta') -
\om}, \g e^{i(\tau + \ta') - \om}, \d e^{i(\ta' + \tau) - \om},
e^{2i\ta'};q)_\ty}\cr
&\cdot \sum^{N-1}_{j=0} {(\b\g e^{-i(\tau + \ta')}/c, qe^{-
i\tau}/\a\d, qe^{-i\ta'}/c;q)_j (qe^{2\om - 2i (\tau +
\ta')};q)_{2j}\o
(q, qe^{\om - i(\tau + \ta')}/\a, qe^{\om - i(\tau + \ta')}/\b,
qe^{\om - i (\tau + \ta')}/\g, qe^{\om -i(\tau + \ta')}/\d;q)_j}\cr
&\cdot \left(- {c \o \b\g} e^{i(\tau + \ta')} \right)^j
q^{-{j\choose 2}}
\ _4\phi_3\left[\matrix{
q^{-j}, ae^{-i\ta'}, be^{-i\ta'}, de^{-i\ta'}\cr
qe^{-2i\ta'}, \b\g e^{-i(\tau + \ta')}/c, \a\d e^{i\tau} q^{-j}\cr}
;q,q\right]\cr
&\cdot P^{\lm, \mu}_{e^{i\tau}} (\cos\ta, \cos (\tau + \ta' + i\om
+ ij \log q))\cr
&+ \ {\rm idem}\ (\ta'; - \ta')\cr
&= {\d (x-x')\o \rho^\lm (x)},\cr
}
$$
where, as before, $x = \cos\ta$, $\g/c = e^\om$, $q^{1-N} < \g/c
< q^{-N}$ and $\a\g = ac$, $\b\d = bd$. When $N=1$, the formula
above reduces to (11.6).
\bigskip
\noindent
{\bf 13. General nonsymmetric $q$-Fourier transformation and its
inversion formula}.\quad From (3.7), (3.8) and the orthogonality
property (2.4) of the Askey--Wilson polynomials one can write
$$t^m \, r^\lm_m (x) = \int^1_{-1} P^{\lm, \mu}_t (x,y) \, r^\mu_m
(y)\,
\rho^\mu (y) \, dy,\leqno (13.1)
$$
where $|t| < 1$ and $r^\lm_m (x) = h^\lm_m \, p^\lm_m (x)$,
$r^\mu_m(y) = h^\mu_m \, p^\mu_m(y)$. To consider the analytic
continuation of
(13.1) to the unit circle $|t| = 1$, when $t = re^{i\tau} =
e^{i\tau - \kappa}$, $\kappa \to 0^+$ and $\g > c$, let us rewrite
this formula in the contour integral form (Fig.~2):
$$t^m \, r^\lm_m (x) = {\log q^{-1}\o i} \, h^\lm_0
\int_{C^-} {L^{\lm,\mu}_t (z, s) \, r^\mu_m \bg{{1\o 2} (q^s + q^{-
s})} \bg{q^{2s}, q^{-2s};q}_\ty ds\o
(\a q^s, \a q^{-s}, \b q^s, \b q^{-s}, \g q^s, \g q^{-s}, \d q^s,
\d
q^{-s};q)_\ty},\leqno (13.2)$$
where $x = \cos\ta, q^z = \eit$ and $y = \cos\phi, q^s = \eip$. The
singular term of $L^{\lm,\mu}_t (z,s) = K^{\lm,\mu}_t (x,y)$ is the
$k = 0$ term of $K^{(3)}_t (x,y)$ in (3.12), and the poles are
defined
by (10.3). Therefore, by analytic continuation, from (13.2) we
obtain
$$e^{i\tau m} \, r^\lm_m (x) = \int^1_{-1}
P^{\lm,\mu}_{e^{i\tau}}(x,y)\,
r^\mu_m (y) \, \rho^\mu (y) \, dy,\leqno (13.3)$$
provided that $\g > c$.
\medskip
To consider the analytic continuation of (13.1) with $\lm$ and
$\mu$ interchanged, $t$ replaced
by $t' = re^{-i\tau} = e^{-i\tau - \kappa}$ and $\kappa \to 0^+$,
let us rewrite this formula in terms of a contour integral
(Fig.~3):
$$(t')^m \, r^\mu_m(y) = {\log q^{-1}\o i} \, h^\mu_0
\int_{C^-} {L^{\mu,\lm}_{t'} (s, z) \, r^\lm_m \bg{{1\o 2} (q^z +
q^{-
z})} \bg{q^{2z}, q^{-2z};q}_\ty dz\o
(aq^z, aq^{-z}, bq^z, bq^{-z}, cq^z, cq^{-z}, dq^z,
dq^{-z};q)_\ty}.
\leqno (13.4)$$
The essential poles of the integrand are now at
$$\eqalign{
z_1 &= {\om-\kappa - i (\tau + \phi)\o
\log q^{-1}}, \quad z_1-1, \ldots, \cr
z_2 &= {\om-\kappa - i (\tau - \phi)\o
\log q^{-1}}, \quad z_2-1, \ldots\ . \cr}\leqno (13.5)
$$
In the limit $\kappa \to 0^+$ the first poles, $z_1$ and $z_2$,
will move from the left half plane to the right half plane, thus
invalidating the formula (13.4). Therefore, in analogy with
section~10,
we need to replace our contour $C^-$ by a contour $C^-_\lm$
with certain indentations (see Fig.~3). This gives the analytic
continuation of (13.4), and for the kernel $P$ we can write
$$e^{-i\tau m} \, r^\mu_m (y) = \int_{\G_\lm}
P^{\mu,\lm}_{e^{-i\tau}}(y,x)
\, r^\lm_m (x) \, \rho^\lm (x) \, dx.\leqno (13.6)$$
Here $C^-_\lm \to \G_\lm$ when $x = {1\o 2} (q^z + q^{-z})$.
\medskip
Now we can define a general $q$-Fourier transformation by
$$f(x) = \int^1_{-1} P^{\lm,\mu}_{e^{i\tau}} (x,y) \, g(y) \,
\rho^\mu(y) \, dy :=
F_q[g](x),\leqno (13.7)$$
and may guess that its inversion formula has the form
$$
g(y) = \int_{\G_\lm} P^{\mu,\lm}_{e^{-i\tau}} (y,x) \,
f(x) \, \rho^\lm (x) \, dx := F_q^{-1}[f](y). \leqno (13.8)
$$
From this point of view, eqs. (13.3) and (13.6) are just mappings,
up to a phase factor, of two systems of Askey--Wilson polynomials,
$\{ p^\lm_m (x)\}$ and $\{p^\mu_m (y)\}$, corresponding to two
different sets of parameters $\lm = (a, b, c, d)$ and $\mu = (\a,
\b,\g, \d)$ with the conditions $\a\g = ac$ and $\b\d = bd$.
\medskip
To prove the inversion formula (13.8), let us interchange the order
of integration in the double integral
$$\eqalignno{
&\int_{\G_\mu} dy \, \rho^\mu (y) \, P^{\lm,\mu}_{e^{i\tau}} (x,y)
\int_{\G_{\lm'}} P^{\mu,\lm'}_{e^{-i\tau}} (y, x') \,
f(x') \, \rho^\lm (x') \, dx'\cr
& \qquad = \int_{\G_{\lm'}} dx' \, \rho^{\lm'}(x') \, f(x')
\int_{\G_\mu}
P^{\lm,\mu}_{e^{i\tau}} (x,y) \, P^{\mu, \lm'}_{e^{-i\tau}} (y, x')
\,
\rho^\mu (y) \, dy\cr
}
$$
and apply the multiplication law (7.1) for the kernels to get
$$\leqalignno{
&\int_{\G_\mu} dy \, \rho^\mu (y)\, P^{\lm,\mu}_{e^{i\tau}} (x,y)
\int_{\G_{\lm'}} P^{\mu,\lm'}_{e^{-i\tau}} (y, x') \,
f(x') \,  \rho^{\lm'}(x') \, dx'&(13.9)\cr
& \qquad = \int^1_{-1} P_1^{\lm, \lm'} (x, x') \, f(x') \,
\rho^{\lm'}(x') \, dx'.\cr
}
$$
In the limit $\lm' \to \lm$, under the same conditions as in
(10.1), we obtain the inversion formula (13.8) for bounded
continuous functions.
\medskip
By using the same considerations as in section 12 one can finally
obtain the inversion formula in terms of a real integral with
additional mass points:
$$\leqalignno{&&(13.10)\cr
&g (y) = \int^1_{-1} P^{\mu,\lm}_{e^{-i\tau}} (y,x) \,
f(x) \, \rho^\lm (x) \, dx\cr
&+{(\a\b, \a\d, \b\d, \a\eip, \b\eip, \g\eip, \d\eip, e^{2i(\tau +
\phi) - 2\om};q)_\ty\o
(ab, ad, bd, ae^{i(\tau + \phi) -\om}, be^{i(\tau + \phi) -\om},
ce^{i(\tau + \phi) -\om}, de^{i(\tau + \phi) -\om},
e^{2i\phi};q)_\ty}\cr
&\cdot \sum^{N-1}_{j=0} {(\g e^{i(\tau-\phi)}, qe^{-i\tau}/\a\d,
qe^{-i\phi}/\b;q)_j (qe^{2\om - 2i(\tau + \phi)};q)_{2j}\o
(q, qe^{\om-i (\tau + \phi)}/a, qe^{\om-i(\tau + \phi)}/b, qe^{\om
- i(\tau +\phi)}/c, qe^{\om - i (\tau + \phi)}/d;q)_j}\cr
&\cdot \left( -\g e^{i(\tau+\phi)} \right)^j q^{-{j\choose 2}}\
_4\phi_3
\left[\matrix{
q^{-j}, ae^{i(\tau + \phi)-\om} q^{-j}, be^{i(\tau + \phi) - \om}
q^{-j}, de^{i(\tau + \phi) -\om} q^{-j}\cr
q^{1-2j} e^{2i (\tau + \phi) - 2\om}, \b\g e^{i\phi-\om} q^{-j}/c,
\a\d e^{i\tau} q^{-j}\cr} ;q,q\right]\cr
&\cdot f (\cos (\tau + \phi + i\om + ij \log q))\cr
&+\ {\rm idem}\ (\phi; -\phi),\cr
}
$$
where $y = \cos\phi, \g/c = e^\om$, $q^{1-N} < \g/c < q^{-N}$, and
$\a\g = ac$, $\b\d = bd$.
\bigskip
\noindent
{\bf 14. Special cases of the kernel}.\quad Although the general
form of the kernel $K_t^{\lm, \mu} (x,y)$ given in (3.9)--(3.12) is
quite formidable it does expose all its poles which are all that
are needed for the purposes of this paper. However, in view of its
multiplicative property (7.2) and its continuous orthogonality when
$|t| \to 1$ it may be useful to list some special cases that
correspond to orthogonal polynomials which are lower than the
Askey--Wilson polynomials in the general scheme of classical
orthogonal polynomials.
\bigskip
\noindent
{\bf Case I}.\ {\it Continuous dual $q$-Hahn polynomials}. \
Setting
$d = 0$ and $\delta = 0$, $\a\g = ac$ in (3.8) one can obtain
$$
\leqalignno{&&(14.1) \cr
&K_t^{\lm,\mu}(x,y) \cr
&:= \sn {(ab, ac;q)_n\o (q, bc;q)_n} \, (ta^{-2})^n \, p_n (x;a, b,
c)\,
p_n(y; \a, \b,\g) \cr
&= { (c^2t^2/\g^2, \a ct \eit/\g, \a ct \meit/\g;q)_\ty \o
(\a\b, ac, bc, \a c^2 t^2 \eip/\g^2;q)_\ty } \cr
&\cdot { (\b \eip, bct \meip/\g, \g \meip, c^2t\eip/ \g, \a
t\eip;q)_\ty\o
( cte^{i\ta + i\phi}/ \g, cte^{i\ta -i\phi}/ \g,
cte^{i\phi - i\ta}/ \g, cte^{-i\ta - i\phi}/ \g ;q)_\ty}\cr
&\cdot \tk (-bc)^k \, q^{k\choose 2}\,
{(\a c^2 t^2 \eip/\g^2 ;q)_{2k} \o
(c^2t^2/\g^2;q)_{2k}} \cr
&\cdot
{( t, cte^{i\ta + i\phi}/\g, cte^{i\ta - i\phi}/ \g,
cte^{i\phi -i\ta}/\g, cte^{-i\ta-i\phi}/\g;q)_k \o
( q, c^2t \eip/\g, \a t \eip, bct \meip/\g,
\a ct \eit/\g, \a ct \meit/\g;q)_k } \cr
&\cdot  \sum^\ty_{\ell = 0} \left(\b \eip \right)^\ell \,
{ (\a \meip, bctq^k/\b\g, \a c^2 t^2 q^{2k}\eip/\g^2;q)_\ell \o
(q, bct q^k \meip/\g, c^2t^2 q^{2k}/\g^2;q)_\ell } \cr
&\cdot { ( ctq^k e^{i\ta - i\phi}/\g, ctq^ke^{-i\ta
-i\phi}/\g;q)_\ell \o
( \a ctq^k \eit/\g, \a ctq^k \meit/\g;q)_\ell} \cr
&\cdot \ _8W_7 \left({\a c^2 t^2 \o \g^2} q^{2k+\ell -1}
\eip; {ct \o \g} q^k e^{i\ta + i\phi}, {ct \o \g} q^k e^{i\phi
-i\ta},
\a \eip, {\a t \o \g} q^{k+\ell}, {c^2t \o \g^2} q^{k +\ell};
q, \g \meip \right).\cr
}
$$
Here $\lm = (a, b, c)$ and $\mu = (\a, \b, \g)$, $\a \g = a c$.
It would be nice if this formula would follow from (3.10)--(3.12)
by
simply setting $d= \d = 0$, but it doesn't, unless one goes through
a series of transformations followed by careful use of the limits.
We shall give a separate proof of (14.1) in the Appendix.
\medskip
The polynomials in (14.1) are, of course, the continuous dual
$q$-Hahn polynomials
$$
p_n(x; a, b,c): = \ _3\phi_2\left[\matrix{
q^{-n},\, a \eit, \, a \meit \cr
\cr
ab, \, ac,\cr} ;\,q, \,q \right],\leqno (14.2)
$$
see, for example, [6], [7] . It is easy to see from (14.1) that in
the limit $t\to 1^-$ the kernel reduces to
$$\leqalignno{
K_1^{\lm, \mu} (x,y)
&= \lim_{t\to 1^-} K^{\lm, \mu}_t (x,y)&(14.3)\cr
&= {(\a e^{i\phi}, \a e^{-i\phi}, \b e^{i\phi}, \b e^{-i\phi},
ce^{i\ta}, ce^{-i\ta}, (c/\g)^2;q)_\ty\o
(\a\b, ac, bc, c e^{i\ta + i\phi}/\g, c e^{i\ta - i\phi}/\g,
c e^{i\phi -i\ta}/\g, c e^{-i\ta -i\phi}/\g;q)_\ty},\cr
}
$$
$c\ne \g$, $\a\g = ac$ and $\b\g = bc$, and we can apply the above
method to get an inversion formula for the corresponding
transformation.
\bigskip
\noindent
{\bf Case II}.\ {\it Al-Salam--Chihara polynomials}.\ The
Al-Salam and Chihara polynomials [6], [12] may be defined by simply
setting
$c = 0$ in (14.2), i.e.
$$\leqalignno{
p_n(x; a, b) &:=\ _3\phi_2\left[\matrix{
q^{-n}, \, a \eit, \, a \meit \cr
\cr
ab, \, 0 \cr} ;q, \, q \right]&(14.4) \cr
&= {(b \meit ;q)_n \o (ab;q)_n} \left( a \eit \right)^n \
_2\phi_1 \left[ \matrix{
q^{-n}, \, a \eit \cr
\cr
q^{1-n} \eit/b \cr} ;q,\, q \meit /b \right].\cr
}
$$
So, with $\lm = (a, b)$ and $\mu = (\a, \b)$, $\a\b = ab$, we find
that
$$\leqalignno{ &&(14.5) \cr
K_t^{\lm,\mu} (x,y)
&: = \sn {(ab;q)_n \o (q;q)_n} \, (ta^{-2})^n \, p_n (x; a, b) \,
p_n (y;\a,\b) \cr
&= { (\a^2 t^2 / a^2, b \meit, \a^2 t \eit/a, b t \eit,
\a t \eip, \a t \meip ;q)_\ty \o
(ab, \a^2 t^2 \eit/a, \a te^{i\ta + i\phi}/a, \a t e^{i\ta
-i\phi}/a,
\a te^{i\phi - i\ta}/a, \a t e^{-i\ta - i\phi}/a;q)_\ty } \cr
&\cdot \ _8W_7\left( \a^2 t^2 \eit/aq; t, \a t/\b, a \eit,
\a te^{i\ta + i\phi}/a, \a t e^{i\ta - i\phi}/a;\,q,\, b \meit
\right),
\cr
}
$$
which in the case $t = 1$ becomes
$$K^{\lm, \mu}_1(x, y) = { ( b \eit, b \meit, \a \eip, \a \meip,
\a^2/a^2;q)_\ty \o
( ab, \a e^{i\ta + i\phi}/a, \a e^{i\ta - i\phi}/a,
\a e^{i\phi - i\ta}/a, \a e^{-i\ta - i\phi}/a;q)_\ty}.\leqno (14.6)
$$
With a different normalization, namely,
$$\leqalignno{
p_n (x; a, b)
&:= {(ab;q)_n\o (q;q)_n} \, a^{-n} \ _3\phi_2 \left[ \matrix{
q^{-n}, \, a \eit, \, a \meit\cr
\cr
ab, \, 0 \cr} ;\,q, \,q \right],&(14.7)\cr
}
$$
the corresponding kernel is
$$\leqalignno{
K_t^{\lm,\mu} (x,y)
&: = \sn {(q;q)_n \o (ab;q)_n} \, t^n \, p_n (x; a, b) \,
p_n (y;\a,\b) &(14.8) \cr
&= { (t^2, b \meit, \a t \eit, \b t \eit, a t \eip, a t \meip
;q)_\ty \o
(ab, a t^2 \eit, te^{i\ta + i\phi}, t e^{i\ta -i\phi},
te^{i\phi - i\ta}, t e^{-i\ta - i\phi};q)_\ty } \cr
&\cdot \ _8W_7\left( a t^2 q^{-1} \eit;  \a t/b, \b t/b, a \eit,
te^{i\ta + i\phi}, t e^{i\ta - i\phi};\,q,\, b \meit \right)
\cr
&= { (\b t/a, \a t \eit, \a t \meit, a t \eip, a t \meip ;q)_\ty \o
(\a at, te^{i\ta + i\phi}, t e^{i\ta -i\phi},
te^{i\phi - i\ta}, t e^{-i\ta - i\phi};q)_\ty } \cr
&\cdot \ _8W_7\left( \a atq^{-1};  \a t/b,  a \eit, a \meit,
\a \eip, \a \meip;\,q,\, \b t/a \right)
\cr
}
$$
with $ab = \a\b$. Eqs. (14.5) and (14.8) can be derived as a
limiting case
of (14.1).
\medskip
The special case of the polynomials (14.7) are the so-called
continuous
$q$-Laguerre polynomials,
$$
L_n^{\a}(x|q) = p_n\left(x; q^{(2 \a + 1)/4},\, q^{(2 \a +
3)/4}\right)
\leqno(14.9)
$$
(see [6]), which go to the classical Laguerre polynomials in the
limit
$q\to 1$,
$$
\lim_{q\to1}L_n^{\a} \left( 1-(1-q^{1/2})x \Big| q \right) =
L_n^{\a}(x).
\leqno(14.10)
$$
From this point of view, one can consider (14.8) as a $q$-version
of
the well-known Poisson kernel for the Laguerre polynomials,
$$\leqalignno{
&\sn {n! \o \G(\a + n + 1)}\, t^n \, L_n^{\a}(x)\,
L_n^{\a}(y)&(14.11)\cr
&\quad= (1-t)^{-1} \exp \left( -t\,{x + y \o 1 - t}\right)\,
\left(xyt \right)^{-\a/2}\,
I_{\a} \left( 2\,{ (xyt)^{1/2} \o 1 - t } \right),
\quad |t|<1, \cr
}
$$
where
$$\leqalignno{
I_{\nu}(z)
&={ (z/2)^{\nu} \o \G (\nu + 1) }\,
\ _0F_1\left( -;\, \nu + 1;\, z^2/4 \right) &(14.12) \cr
&={ (z/2)^{\nu} \o \G (\nu + 1) }\, \ e^{-z}\,
\ _1F_1\left( \nu + 1/2;\, 2 \nu + 1;\, 2 z \right). \cr
}
$$
So, one can consider the corresponding transformation with the
kernel (14.8)
as a $q$-version of the Fourier--Bessel transform.
\bigskip
\noindent
{\bf Case III}.\ {\it Continuous big $q$-Hermite polynomials}.\
Setting $b = 0$ in (14.4) one obtains the so-called continuous big
$q$-Hermite polynomials [12]:
$$\leqalignno{
p_n (x;a)
&:=\ _3\phi_2\left[\matrix{
q^{-n}, \, a \eit,\, a \meit \cr
\cr
0,\, 0 \cr} ;\,q,\,q \right] &(14.13)\cr
&= \left(a \eit \right)^n\,
\sum^n_{k=0}\, {(q^{-n}, a \eit;q)_k \o (q;q)_k }\,
\left( -q^n e^{-2i\ta}\right)^k\, q^{-{k\choose 2}}.\cr
}
$$
The nonsymmetric Poisson kernel for these polynomials can be
derived
as a limiting case of (14.5) by [10, (3.2.11)],
$$\leqalignno{
K^{a, \a}_t (x,y)
&:= \sn {(ta^{-2})^n\o (q;q)_n}\, p_n(x;a)\, p_n(y;\a)&(14.14)\cr
&= { \left( \a^2 t^2/a^2, \a t\eip, \a \meip;q \right)_\ty \o
\left( \a te^{i\ta + i\phi}/a, \a te^{i\ta - i\phi}/a,
\a te^{i\phi - i\ta}/a, \a te^{-i\ta - i\phi}/a ;q \right)_\ty}\cr
&\cdot\ _3\phi_2\left[\matrix{
t, \, \a te^{i\ta + i\phi}/a, \, \a te^{i\phi - i\ta}/a \cr
\cr
\a^2 t^2/a^2,\, \a t \eip \cr} ;\,q, \, \a \meip \right].\cr
}
$$
\indent
On the other hand, setting $b = \b = 0$ and $b/\b = \a/a$
in (14.7)--(14.8) one obtains the kernel
$$\leqalignno{
K^{a,\a}_t (x,y)
&:= \sn t^n\, (q;q)_n\, p_n(x;a)\, p_n(y;\a)&(14.15)\cr
&= {(t^2, at \eip, \a \meip;q)_\ty\o
(te^{i\ta + i\phi}, te^{i\ta - i\phi},
te^{i\phi - i\ta}, te^{-i\ta-i\phi};q)_\ty}\cr
&\cdot\ _3\phi_2 \left[\matrix{
at/\a,\, te^{i\ta + i\phi},\, te^{i\phi - i\ta}\cr
\cr
t^2, \, at \eip \cr} ;\,q,\, \a \meip \right].\cr
}
$$
Observe that
$$\lim_{a\to 0} \, e^{in\ta} \, \sn {(q^{-n}, a \eit ;q)_k \o
(q;q)_k} \,
\left(-e^{-2i\ta}\right)^k q^{-{k\choose 2}} = H_n (x|q),\leqno
(14.16)$$
where $H_n (x|q)$ is the continuous $q$-Hermite polynomial defined
in (1.10). See also [1] for the proof of this limit directly from
the first expression in (14.13).
Equation (14.15) then reduces further to
$$\leqalignno{
K^{0,\a}_t (x,y)
&= {(t^2, \a \meip;q)_\ty\o
(te^{i\ta + i\phi}, te^{i\ta - i\phi}, te^{i\phi - i\ta}, te^{-i\ta
- i\phi};q)_\ty}&(14.17)\cr
&\cdot\ _2\phi_1 \left[\matrix{
te^{i\ta + i\phi}, \, te^{i\phi -i\ta}\cr
\cr
t^2\cr} ; \,q, \, \a \meip \right].\cr
}
$$
This $_2\phi_1$ series is related to an addition formula [16] for
the $q$-Bessel functions of Jackson [10, Ex. 1.24]. Denoting
$$J_t (x,y)=\ _2\phi_1\left[\matrix{
te^{i\ta + i\phi}, \, te^{i\ta - i\phi}\cr
\cr
t^2\cr}; q, \a \meit \right],\leqno (14.18)
$$
we find that a special limiting case of (10.6) is
$$\leqalignno{
&\int_\Gamma J_{e^{i\tau}}(x,y) \, J_{e^{i\tau}}^* (x', y) \,
\rho \left(y; e^{i\tau + i\ta}, e^{i\tau - i\ta},
e^{i\ta' - i\tau}, e^{-i\ta'- i\tau} \right)\, dy&(14.19) \cr
&\qquad= \left[ {2\pi \o (q;q)_\ty} \right]^2 \,
{ \delta (x-x') \o \Big| \Big( e^{2i\tau}, e^{2i\ta};q \Big)_\ty
\Big|^2 }
\, (1-x^2)^{1/2},\cr
}
$$
where $\G$ is a limiting symmetric form of the contour $\G_\mu$
in (10.6).
\medskip
We would like to add, in conclusion, that there are many more
interesting special and limiting cases of the Askey--Wilson
polynomials and that the same technique as we have developed in
this paper can be applied to most of them.
\bigskip
\noindent
{\bf 15. Appendix: proof of (14.1)\/}.\quad  The $q$-integral
representation
(4.4) in the case of the dual $q$-Hahn polynomials has the simpler
form
$$\leqalignno{
K_t^{\lm, \mu} (x,y) = B_1(\ta, \phi)
&\int^{q\meit/b}_{q\eit/b}d_qu\,
{(bu\eit, bu\meit;q)_\ty \o
(bau/q, bcu/q;q)_\ty} &(15.1)\cr
&\cdot \int^{q\meip/\g}_{q\eip/\g} d_qv \, {(\g v \eip, \g
v\meip;q)_\ty \o
(\g\a v/q, \g\b v/q;q)_\ty} \cr
&\cdot\, _2\phi_1  \left[ \matrix{
q/u,\, q/v \cr
bc \cr} ;\,q,\,{bc \o q^2}\,uvt \right],\cr
}
$$
where
$$\leqalignno{
B^{-1}_1 (\ta,\phi) = - \  & {q^2 (1-q)^2\o 4b\g} \,(q, q, ac,
\a\b;q)_\ty
&(15.2)\cr
&\cdot h(x;b)\, h(y;\g)\, \rho^\lm (x)\, \rho^\mu (y),\cr
}
$$
and $\lm = ( a, b, c)$ and $\mu = (\a, \b, \g)$, $\a \g = ac$.
The troublesome $\ _6W_5$ function in (4.4) now becomes a $\
_2\phi_1$
series which is first transformed to a $\ _2\phi_2$ series by
$$
\, _2\phi_1  \left[ \matrix{
q/u,\, q/v \cr
bc \cr} ;\,q,\,{bc \o q^2}\,uvt \right]
={ (bctu/q, bctv/q;q)_\ty \o (bc, bcuvt/q^2;q)_\ty }
\cdot\, _2\phi_2  \left[ \matrix{
t, \, bcuvt/q^2 \cr
bctu/q,\, bctv/q \cr} ;\,q,\,bc \right]
\leqno (15.3)
$$
(this formula can be obtained as a consequence of [10, III.2 and
III.4]).
So,
$$
K_t^{\lm, \mu}(x, y) = {B_1 \o (bc;q)_\ty} \,
\tk (-bc)^k \, q^{ {k\choose 2} }\, {(t;q)_k \o (q;q)_k} \, V_k,
\leqno (15.4)
$$
where
$$\leqalignno{
V_k = &  \int^{q\meit/b}_{q\eit/b} d_q u \,
{(bu\eit, bu\meit, bctuq^{k-1};q)_\ty \o
(bau/q, bcu/q;q)_\ty} &(15.5)\cr
&\cdot \int^{q\meip/\g}_{q\eip/\g}d_q v \,
{(\g v \eip, \g v\meip, bctvq^{k-1};q)_\ty \o
(\g\a v/q, \g\b v/q, bcuvtq^{k-2};q)_\ty}.\cr
}
$$
The limiting form of the $q$-integral (6.1) now is
$$\leqalignno{
&\int^{q\meip/\g}_{q\eip/\g} {(\g v\eip, \g v\meip, bct
vq^{k-1};q)_\ty \o
(\g\a v/q, \g\b v/q, bcuvtq^{k-2};q)_\ty} \, d_qv &(15.6)\cr
&= {q(1-q)\o 2i\g} \, (q, \a \b ;q)_\ty \, h(y;\g) \, \rho^\mu(y)
\,{(bctq^k\meip/\g;q)_\ty\o (bcutq^{k-1}\meip/\g ;q)_\ty}\cr
&\cdot\, _3\phi_2  \left[ \matrix{
\a \meip, \, \b \meip, \, q/u \cr
\a\b, \, bctq^k\meip/\g \cr} ;\,q,\,{bcut \o \g}\,q^{k-1}\eip
\right] \cr
&= {q(1-q)\o 2i\g} \, (q;q)_\ty \, h(y;\g) \, \rho^\mu(y)
\,{(\a\eip, \a\meip, bctq^k/\g\b, bc \b utq^{k-1}/\g;q)_\ty \o
(\a/\b, bcutq^{k-1}\eip/\g, bcutq^{k-1}\meip/\g ;q)_\ty}\cr
&\cdot\, _3\phi_2  \left[ \matrix{
bctq^k/\a\g, \,\b\eip, \, \b\meip \cr
q\b/\a, \, bc \b utq^{k-1}/\g \cr} ;\,q,\, q \right] +
\ {\rm idem} (\a; \b), \cr
}
$$
see also [10, Ex. 3.6] for the last transformation. Therefore,
$$\leqalignno{
K_t^{\lm, \mu}(x,y)
&={(\a\eip, \a\meip, bct/\g\b;q)_\ty \o
(ab, bc, \a\b, \a/\b;q)_\ty} &(15.7) \cr
&\cdot {(ct\eit, c\b t\eit/\g, a\meit, c\meit ;q)_\ty \o
(cte^{i\ta + i\phi}/\g, cte^{i\ta - i\phi}/\g, e^{-2i\ta} ;q)_\ty
} \cr
&\cdot \tk
{(t,cte^{i\ta + i\phi}/\g, cte^{i\ta - i\phi}/\g ;q)_k \o
(q, bct/\g\b, ct\eit, c\b t\eit/\g ;q)_k}\, q^{{k\choose 2}}\,
(-bc)^k \cr
&\cdot \sum^\ty_{\ell = 0} {(bctq^k/\a\g, \b \eip, \b \meip;q)_\ell
\o
(q, q\b/\a, c\b tq^k\eit/\g;q)_\ell}\, q^\ell \cr
&\cdot \, _4\phi_3  \left[ \matrix{
ctq^ke^{i\ta + i\phi}/\g,\, ctq^ke^{i\ta - i\phi}/\g,\, a\eit,\,
c\eit \cr
qe^{2i\ta}, \, ctq^k\eit,\, c \b tq^{k + \ell}\eit/\g \cr } ;
\,q,\,q \right] \cr
& + \ {\rm idem} (\ta; -\ta) + \ {\rm idem} (\a; \b). \cr
}
$$
Unfortunately, the $\ _4\phi_3$ series appearing in (15.7) are not
balanced and cannot be transformed in a simple way.
By a change of order of summation in the last two series one gets
$$\leqalignno{
K_t^{\lm, \mu}(x,y)
&={(ct\eit, a\meit, c\meit;q)_\ty \o
(ac, bc, \a\b,
cte^{i\ta + i\phi}/\g, cte^{i\ta - i\phi}/\g, e^{-2i\ta} ;q)_\ty}
&(15.8) \cr
&\cdot \tk
{(t,cte^{i\ta + i\phi}/\g, cte^{i\ta - i\phi}/\g ;q)_k \o
(q, ct\eit;q)_k}\,
q^{{k\choose 2}}\, (-bc)^k \cr
&\cdot \sum^\ty_{\ell = 0}
{(ctq^ke^{i\ta + i\phi}/\g, ctq^ke^{i\ta - i\phi}/\g, a\eit,\,
c\eit; q)_\ell
\o (q, qe^{2i\ta}, ctq^k\eit; q)_\ell}\, q^\ell \cr
&\cdot \Bigg\{ {(\a\eip, \a\meip, bctq^k/\b\g, \b
ctq^{k+\ell}\eit/\g;q)_\ty
\o (\a/\b;q)_\ty} \cr
&\cdot \, _3\phi_2  \left[ \matrix{
bctq^k/\a\g,\, \b \eip,\, \b\meip \cr
q\b/\a, \, \b ctq^{k+\ell}\eit/\g \cr } ;\,q,\,q \right]
+ \ {\rm idem} (\a; \b) \Bigg\} \cr
& + \ {\rm idem} (\ta; -\ta). \cr
}
$$
We can now transform the pairs of $\ _3\phi_2$ series above by [10,
III.34],
thus getting
$$\leqalignno{
K_t^{\lm, \mu}(x,y)
&={(ct\eit, a\meit, c\meit, \a ct\eit/\g, \b\eip,
bct\meip/\g;q)_\ty \o
(ac, bc, \a\b,
cte^{i\ta + i\phi}/\g, cte^{i\ta - i\phi}/\g, e^{-2i\ta} ;q)_\ty}
&(15.9) \cr
&\cdot \tk
{(t,cte^{i\ta + i\phi}/\g, cte^{i\ta - i\phi}/\g ;q)_k \o
(q, ct \eit, \a ct \eit/\g, bct \meip/\g ;q)_k}\,
q^{{k\choose 2}}\, (-bc)^k \cr
&\cdot \sum^\ty_{\ell = 0}
{(ctq^ke^{i\ta + i\phi}/\g, ctq^ke^{i\ta - i\phi}/\g, a\eit,\,
c\eit; q)_\ell
\o (q, qe^{2i\ta}, ctq^k\eit, \a ctq^k\eit/\g; q)_\ell}\, q^\ell
\cr
&\cdot \, _3\phi_2  \left[ \matrix{
\a\meip,\, ctq^{k+\ell}e^{i\ta - i\phi}/\g,\, bctq^k/\b\g \cr
\a ctq^{k+\ell}\eit/\g,\, bctq^k\meip/\g  \cr } ;\,q,\,\b\eip
\right] \cr
& + \ {\rm idem} (\ta; -\ta). \cr
}
$$
We would like to point out that it is possible to obtain (15.9) as
a limiting case of (3.12) by using [10, (3.2.11)].
\medskip
Changing the order of summation one can  rewrite the sum over
$\ell$
in (15.9) in the following form
$$\eqalign{
&\cdot \sum^\ty_{\ell = 0}
{(\a\meip, bctq^k/\b\g, ctq^ke^{i\ta - i\phi}/\g ;q)_\ell \o
(q, bctq^k\meip/\g, \a ctq^k\eit/\g;q)_\ell} \,
\left(\b\eip\right)^\ell \cr
&\cdot\, _4\phi_3\left[\matrix{
ctq^ke^{i\ta + i\phi}/\g, \, ctq^{k + \ell}e^{i\ta - i\phi}/\g, \,
a\eit, \, c\eit \cr
qe^{2i\ta}, \, ctq^k\eit,\, \a ctq^{k+\ell}\eit/\g \cr}; \,q, \,q
\right],\cr
}
$$
where the $\,_4\phi_3$ function is now balanced. The two
$\,_4\phi_3$
series that arise from (15.9) in this way can be combined into a
single
$\,_8\phi_7$ series by [10, III.36]. Substituting this result in
(15.9)
with the aid of [10, III.23] we finally get (14.1).
\bigskip
\noindent
{\bf Acknowledgements}.\quad We would like to thank Prof. G. Gasper
for his valuable comments on the first draft of this paper. One of
us (S.K. Suslov) gratefully thanks the Department of Mathematics
and Statistics of Carleton University, Centre de Rechearches
Math\'ematique Universit\'e de Montr\'eal and the Department of
Mathematics of University of Wisconsin--Madison for their
hospitality.
\bigskip
%\vfill
%\eject
\noindent
{\bf References}
\medskip
\item{1.} Askey, R. (1994). {\it Difference operators and
orthogonal
polynomials\/}, Rendiconti di Matematica, to appear.
\medskip
\item{2.} Askey, R., Atakishiyev, N.M. and Suslov, S.K. (1993).
{\it An analog of the Fourier transformation for the $q$-harmonic
oscillator\/}, Preprint IAE-5611/1, Kurchatov Institute, Moscow;
Symmetries in Science VI (B. Gruber, ed.), Plenum Press, New York,
pp.~57--63.
\medskip
\item{3.} Askey, R. and Ismail, M.E.H. (1983). {\it A
generalization of ultraspherical polynomials\/}, Studies in Pure
Mathematics (P. Erd\"os, ed.), Birkh\"auser, Boston, Mass., pp.
55--78.
\medskip
\item{4.} Askey, R. and Suslov, S.K. (1993). {\it The $q$-harmonic
oscillator and an analogue of the Charlier polynomials\/},
J. Phys. A: Math. Gen., {\bf 26} \#15, L693--L698.
\medskip
\item{5.} Askey, R. and Suslov, S.K. (1993). {\it The $q$-harmonic
oscillator and the Al-Salam and Carlitz polynomials\/}, Lett. Math.
Phys., {\bf 29} \#2, 123--132.
\medskip
\item{6.} Askey, R. and Wilson, J.A. (1985). {\it Some basic
hypergeometric polynomials that generalize Jacobi polynomials\/},
Memoirs Amer. Math. Soc. \#319.
\medskip
\item{7.} Atakishiyev, N.M., Rahman, M. and Suslov, S.K. (1994).
{\it On the classical orthogonal polynomials\/}, Constructive
Approximation, to appear.
\medskip
\item{8.} Atatishiyev, N.M. and Suslov, S.K. (1990).
{\it Difference analogs of the harmonic oscillator\/}, Theoretical
and
Mathematical Physics, {\bf 85} \#1, 1055--1062.
\medskip
\item{9.} Biedenharn, L.C. (1989). {\it The quantum group $SU_q(2)$
and a $q$-analogue of the boson operators\/}, J. Phys. A: Math.
Gen.,
{\bf 22} \#18, L873--L878.
\medskip
\item{10.} Gasper, G. and Rahman, M. (1990). {\it Basic
Hypergeometric Series\/}, Cambridge University Press.
\medskip
\item{11.} Gasper, G. and Rahman, M. (1986). {\it Positivity of the
Poisson kernel for the continuous $q$-Jacobi polynomials and some
quadratic transformation formulas for basic hypergeometric
series\/},
SIAM J. Math., Anal., {\bf 17}, 970--999.
\medskip
\item{12.} Koekoek, R. and Swartow, R.F. (1994). {\it The
Askey--scheme of hypergeometric orthogonal polynomials and its
$q$-analogue\/}, Report 94--05, Delft University of Technology.
\medskip
\item{13.} Landau, L.D. and Lifschitz, E.M. (1965). {\it Quantum
Mechanics: Non-relativistic Theory\/}, 2nd ed. (in English),
Addison-Wesley.
\medskip
\item{14.} Macfarlane, A.J. (1989). {\it On $q$-analogues of the
quantum harmonic oscillator and the quantum group $SU(2)_q$\/},
J. Phys. A: Math. Gen., {\bf 22} \#21, 4581--4588.
\medskip
\item{15.} Rahman, M. (1981). {\it The linearization of the product
of continuous $q$-Jacobi polynomials\/}, Canad. J. Math. {\bf 33},
255--284.
\medskip
\item{16.} Rahman, M. (1988). {\it An addition theorem and some
product formulas for $q$-Bessel function\/}, Can. J. Math. {\bf
40},
1203--1221.
\medskip
\item{17.} Rahman, M. and Suslov, S.K. (1993). {\it Singular
analogue of the Fourier transformation for the Askey--Wilson
polynomials\/}, Preprint CRM--1915, October 1993, Centre de
Recherches
Math\'ematique Universit\'e de Montr\'eal, Montr\'eal.
\medskip
\item{18.} Rahman, M. and Verma, A. (1991). {\it Positivity of the
Poisson kernel for the Askey--Wilson polynomials\/}, Indian J.
Math.
{\bf 33} \#3, 287--306.
\medskip
\item{19.} Szeg\"o, G. (1975). {\it Orthogonal Polynomials\/}, 4th
ed., Amer. Math. Soc. Coll. Publ., Vol. 23.
\medskip
\item{20.} Titchmarch, E.C. (1948). {\it Introduction to the Theory
of Fourier Integrals\/}, 2nd ed., Oxford University Press.
\medskip
\item{21.} Wiener, N. (1933). {\it The Fourier Integral and Certain
of Its Applications\/}, Cambridge University Press; Dover edition
published in 1948.
\bye